\documentclass[12pt]{amsart}
\usepackage[margin=1.1in]{geometry}
\usepackage{subcaption, graphicx}
\usepackage{xcolor}
\usepackage{amssymb, amsthm, amsaddr}
\usepackage{tikz}
\usepackage{placeins, url}

\newtheorem{remark}{Remark}
\theoremstyle{definition}
\newtheorem{definition}{Definition}

\definecolor{mycolor1}{rgb}{0.12, 0.47, 0.71}
\definecolor{mycolor2}{rgb}{1.00, 0.73, 0.47}
\definecolor{mycolor3}{rgb}{0.84, 0.15, 0.16}
\definecolor{mycolor4}{rgb}{0.55, 0.34, 0.29}
\definecolor{mycolor5}{rgb}{0.97, 0.71, 0.82}
\definecolor{mycolor6}{rgb}{0.74, 0.74, 0.13}
\definecolor{mycolor7}{rgb}{0.62, 0.85, 0.90}

\title{Topological analysis of U.S. city demographics}

\author{Jakini A. Kauba}
\address{School of Mathematical and Statistical Sciences, Clemson University}
\author{Thomas Weighill}
\address{Department of Mathematics and Statistics, University of North Carolina at Greensboro}

\begin{document}

\begin{abstract}
    We apply persistent homology, the main method in topological data analysis, to the study of demographic data. Persistence diagrams efficiently summarize information about clusters or peaks in a region's demographic data. To illustrate how persistence diagrams can be used for exploratory analysis, we undertake a study of the 100 largest U.S.~cities and their Black and Hispanic populations. We use our method to find clusters in individual cities, determine which cities are outliers and why, measure and describe change in demographic patterns over time, and roughly categorize cities into distinct groups based on the topology of their demographics. Along the way, we highlight the advantages and disadvantages of persistence diagrams as a tool for analyzing geospatial data.
\end{abstract}

\maketitle

\section{Introduction}

The U.S. Census Bureau publishes a wealth of data about the demographic makeup of the United States resulting from the decennial census or from the American Community Survey (ACS). These valuable data allow geographers, sociologists and others to gain insights into how populations are arranged in space, and how those patterns change over time. To do so, the large amount of raw data needs to be quantified via a concise summary. As an example, the study of racial segregation using census data is well established in social science research (see e.g.~\cite{wilson2011visualizing, wright2014patterns}), and relies substantially on measuring segregation using just a single number~\cite{massey1988dimensions}. Notably, not all measures of segregation take spatial arrangement into account -- the Dissimilarity Index~\cite{james1985measures}, for example, takes as input only the list of population values per geographic unit (e.g. census tract) and therefore ignores the inherent spatiality of the data~\cite{brown2006spatial}. 

In this paper we demonstrate a method, based on topological data analysis (TDA), which constructs a summary (the persistence diagram) that quantifies the demographic patterns in a region while also linking the information (i.e.~points) in the persistence diagram to specific locations for further analysis. TDA is an area of mathematics consisting of data analysis methods which rely on theory from algebraic topology to discover the ``shape'' of data, which can remain hidden to conventional machine learning or data science methods. TDA has found applications in a range of areas including medical research, image classification, and data vizualization. It has taken some time for TDA to be used in geospatial applications, but recently TDA methods have been deployed to study gerrymandering \cite{duchin2021homological}, voting patterns~\cite{feng2021persistent} and epidemiological data~\cite{hickok2022analysis}. As for applications to demographic data, there is the recent application of the Mapper algorithm (introduced in \cite{nicolau2011topology}) to Mecklenburg County in \cite{husain2022mappering}. 

In this paper, we take the first step in applying persistent homology to the analysis of demographic data. We use persistent homology to detect and visualize the presence and magnitude of ``peaks'' in the population share of a particular minority group (in our case, either Black or Hispanic). Such peaks correspond to high local concentration of a particular population. The number of peaks and the values near and between them tells us a lot about how a particular population is distributed. The output of the persistent homology method on a particular region is a persistence diagram which summarizes this information. 

To demonstrate how persistence diagrams can help us analyze demographic trends, we undertake a study of the 100 largest U.S.~cities by population using Black and Hispanic population data. Through this study, we demonstrate the usefulness of persistent homology as an exploratory analysis tool. In particular, we use persistent homology to 
\begin{itemize}
    \item[(1)]  find clusters in individual cities (Section \ref{sec:single}),
    \item[(2)] determine which cities are outliers and why (Section \ref{sec:scatter}),
    \item[(3)] measure and describe change over time (Section \ref{sec:change}), and
    \item[(4)] roughly categorize cities into distinct groups (Section \ref{sec:clustering}).
\end{itemize}
The applications (1)--(4) are partly inspired by the applications for redistricting in \cite{duchin2021homological} which follow a similar framework, replacing ``cities'' with ``redistricting plans''.

The outline of the paper is as follows. We begin in Section \ref{sec:background} with an overview of persistent homology, focusing only on what is needed to define the $0^{th}$-dimensional persistence diagram used in our method. In Section \ref{sec:method} we outline our method for turning a city into a persistence diagram. Sections  \ref{sec:single}--\ref{sec:clustering} contain applications (1)--(4) above. Section \ref{sec:drawback} highlights a major drawback of the method, the modifiable areal unit problem (MAUP), after which we conclude our work in the final section.

\section{Background} \label{sec:background}

\subsection{Persistent homology}

For a general mathematics reference for persistent homology we refer the reader to \cite{otter2017roadmap}. For our purposes, we will focus on $0^{th}$-dimensional homology which allows the less technical description below. The basic input to persistent homology is a sequence of spaces 
$$
X_{t_1} \subseteq X_{t_2} \subseteq \cdots \subseteq X_{t_n}
$$
For us, the main example of a space will be a graph, so we can think of an increasing sequence of graphs. For the $0^{th}$ dimensional case, we are interested in counting the connected components in each $X_{t_1}$, and tracking the birth and death of these components. If a new connected component appears in $X_{t_i}$ that was not in $X_{t_{i-1}}$ we say that the component was \emph{born at $t_i$}. If this connected component is joined to a component  with lower birth time (i.e.~an ``older'' component) at $X_{t_j}$, we say that the component \emph{dies at $t_j$}. The connected components remaining in the final space $X_{t_n}$ are traditionally assigned a death time of $\infty$. Once we have collected the (birth time, death time) pairs, we plot each pair as a point in the plane to produce a \emph{persistence diagram}. The \emph{persistence} of a point $(b,d)$ is defined as $d-b$, and high persistence points are usually assumed to be more important or necessary of further investigation.

\subsection{Sublevel set filtrations}

There are many ways to construct a sequence of spaces from data. For example, in the case of point-cloud data, the Vietoris-Rips complex is a common technique. We will employ a different method, namely \emph{sublevel set filtrations}. This technique is common in applications to image classification \cite{chung2021persistent} and has also previously been used on geospatial data \cite{feng2021persistent}. 

Given a graph $G = (V, E)$ and a function on the vertices $f: (V, E) \to \mathbb{R}$, we define $X_{t}$ to be the full subgraph of $G$ whose vertices satisfy $f(v) \leq t$. By \emph{full subgraph}, we mean that all edges in $G$ are added to $X_t$ so long as their endpoints are in $X_t$. Clearly if $s < t$, then $X_s \subseteq X_t$, so if the image of $f$ is the ordered set of values $t_1 \leq t_2 \leq \cdots \leq t_n$ then we have our sequence 
$$
X_{t_1} \subseteq X_{t_2} \subseteq \cdots \subseteq X_{t_n}
$$
from which we can compute a persistence diagram. Points in the persistence diagrams correspond to local minima -- vertices with values lower than all their neighbors.

\subsection{Dual graphs}\label{sec:dualgraph}

\begin{figure}
    \centering
    \begin{subfigure}{0.4\textwidth}
        \includegraphics[width=\textwidth]{ 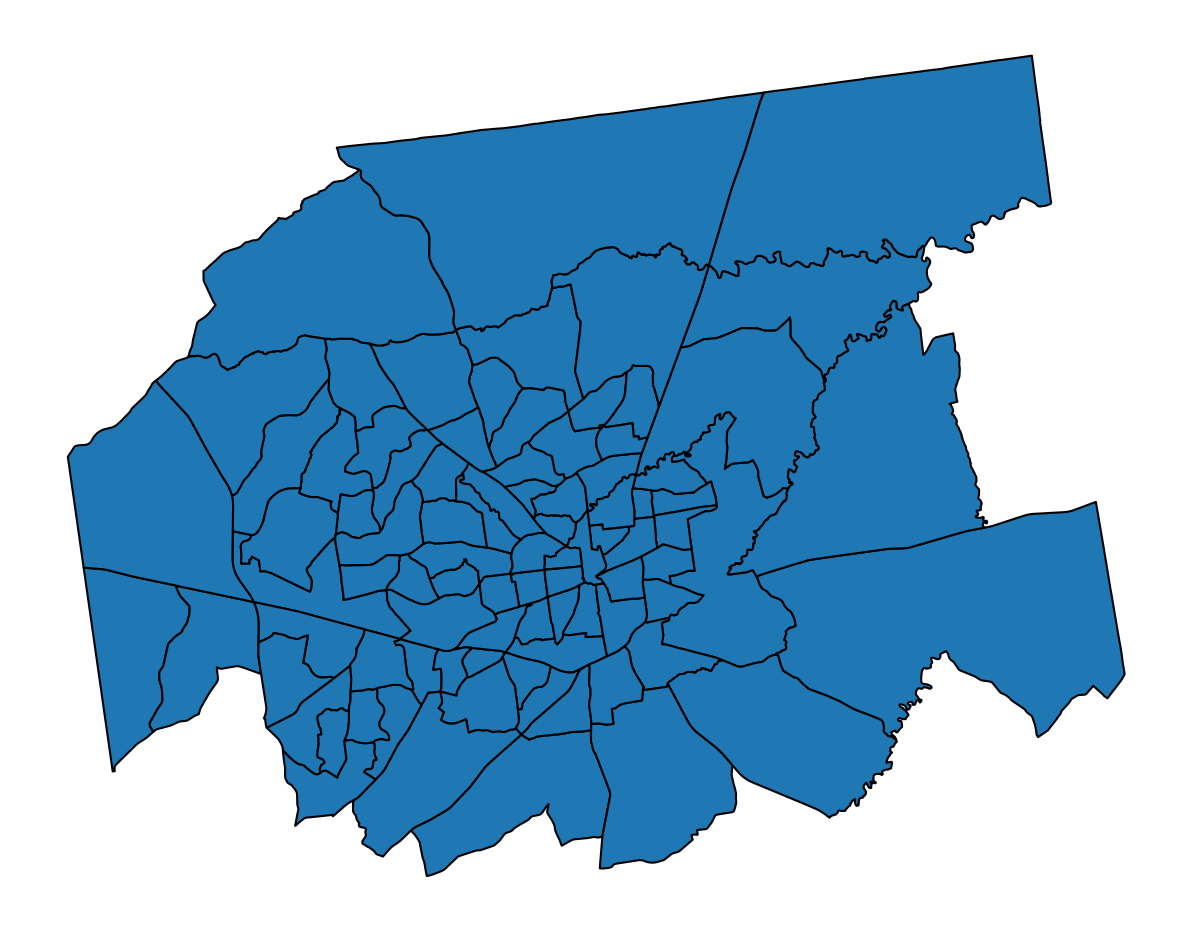}
    \end{subfigure}%
    \begin{subfigure}{0.4\textwidth}
        \includegraphics[width=\textwidth]{ 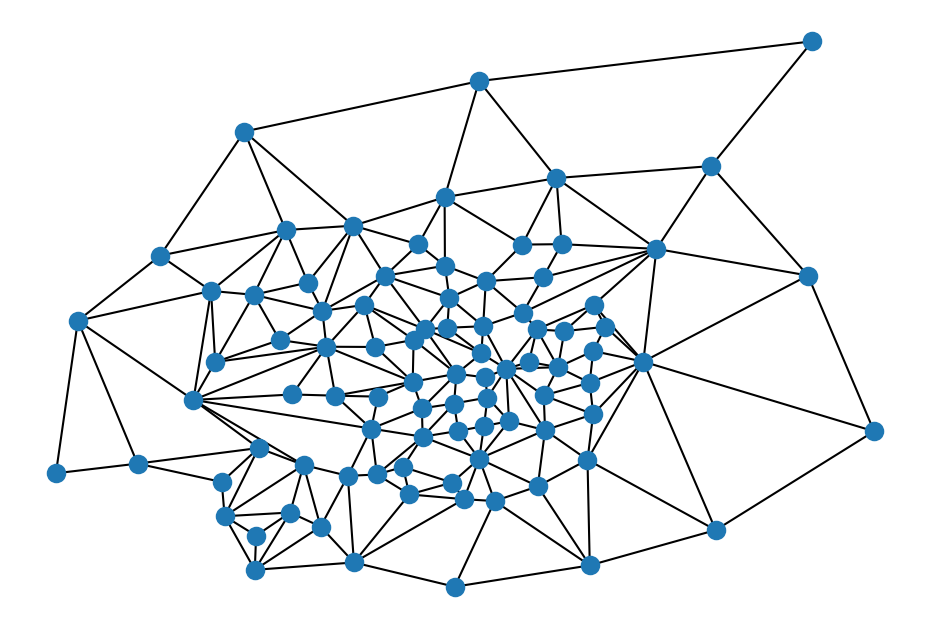}
    \end{subfigure}
    \caption{The census tracts of Greensboro NC and the corresponding dual graph.}
    \label{fig:greensboronc}
\end{figure}

Given geographic units and corresponding values (e.g.~population shares for a minority group), we define the \emph{dual graph} as the graph whose vertices represent geographic units and where there is an edge between every pair of adjacent units. Importantly, the dual graph encodes connectedness of regions: a region is connected if and only if the corresponding subgraph is connected. For this reason dual graphs are particularly useful in redistricting applications where connectedness is a key criterion for districts \cite{duchinAAAS, deford2019recombination}. In Figure~\ref{fig:greensboronc} is an example of the census tracts for the city of Greensboro and the corresponding dual graph. 

\subsection{Wasserstein distance and total persistence}
In order to compare persistence diagrams to one another, we require a notion of distance between persistence diagrams. The established choice in TDA is the family of Wasserstein distances. Given two persistence diagrams $P, Q$, each of which is a finite subset of points in the plane, possibly with a value of infinity in one coordinate, we define a \emph{partial bijection} to be a bijection $\pi: A \to B$, where $A \subseteq P$ and $B \subseteq Q$. In order words, a partial bijection is a matching between some of the points in $P$ and some of the points in $Q$.

\begin{definition}[Wasserstein Distance] Let $P$ and $Q$ be persistence diagrams. We define the \emph{$p$-Wasserstein distance} as follows:
\[
W_p(P,Q)^p = \min_{\pi:A \rightarrow B} \sum_{\alpha \in A}d_p(\alpha, \pi(\alpha))^p + \sum_{\alpha \in P, \alpha \notin A}d_p(\alpha, \Delta)^p + \sum_{\beta \in Q, \beta \notin B}d_p(\beta,\Delta)^p
\]
where $\Delta = \{ (x,x) | x \in \mathbb{R} \}$ and $\pi$ ranges over all partial bijections.
\end{definition}

In the above definition, the distance $d_p$ between points is defined as
$$
d_p((b,d), (b', d')) = \left( |b-b'|^p + |d-d'|^p \right)^{1/p} 
$$
For $p = \infty$, one makes the usual replacement of sums by maxima to obtain
$$
d_\infty((b,d), (b', d')) = \max ( |b-b'|, |d-d'| )
$$
We now mention the $\infty$-Wasserstein distance, also called the bottleneck distance.

\begin{definition}[Bottleneck Distance]\label{def:bnd} Let $P$ and $Q$ be persistence diagrams. We define the \emph{$\infty$-Wasserstein distance} or \emph{bottleneck distance} as follows:
\[
W_\infty(P,Q) = \min_{\pi:A \rightarrow B} \max\left( \max_{\alpha \in A}d_\infty(\alpha, \pi(\alpha)), \max_{\alpha \in P, \alpha \notin A}d_\infty(\alpha, \Delta), \max_{\beta \in Q, \beta \notin B}d_\infty(\beta,\Delta) \right)
\]
where $\Delta = \{ (x,x) | x \in \mathbb{R} \}$ and $\pi$ ranges over all partial bijections.
\end{definition}

Persistence diagrams are complicated objects themselves, and so we may sometimes want to summarize a diagram with a single number. Such numbers are sometimes called persistence statistics~\cite{chung2021persistent}. In this paper, we make use of total persistence as our main statistic.

\begin{definition}[Total persistence]\label{def:tp}
For a persistence diagram $P$, we define the \emph{total persistence} to be the sum of all lifespans:
$$
\mathsf{TP}(P) = \sum_{(b,d) \in P} (d - b)
$$
\end{definition}

In order for $\mathsf{TP}(P)$ to be finite, it is necessary to replace any $\infty$ values with a finite number, such as the maximum possible death time. In this paper we use the value of $1$ for this purpose, since all (finite) death values naturally fall in the interval $[0,1]$.

\section{Methods}\label{sec:method}

\subsection{Data}

As our dataset, we consider the 100 largest U.S.~cities using boundaries obtained from the CDC's 500 Cities project~\cite{cdc}. Census tracts for 2010 and 2020 were obtained from NHGIS~\cite{NHGIS}. Since tracts need not nest perfectly inside city boundaries, we take, for each city, all those tracts which intersect the region inside the city boundaries. To estimate Black population percentage, we assign to each tract the number of persons who indicated their race as Black on the Census divided by the total population of that tract. Similarly, we estimate Hispanic population percentage as the population who indicated their ethnicity as Hispanic divided by the total population. 

\begin{remark}
Some tracts have small populations, making estimates of the Black and Hispanic percentage unreliable or not significant. For example, there was at least one case in which a tract had a total population of 1 and a Black population of 1 (i.e.~100\%). Since our method detects peaks in the Black or Hispanic population percentage, this kind of outlier can create misleading results. We make the somewhat arbitrary choice to treat any tract with total population of less than 10 as a missing data point, and assign to that tract the highest percentage among its neighbors, thus guaranteeing that these empty tracts cannot either form or separate clusters.
\end{remark}

\subsection{Constructing a persistence diagram from city data} \label{sec:constructing}

We now describe how to convert the demographic data for a region into a persistence diagram. Given a set of geographic units, we form their dual graph as defined in Section \ref{sec:dualgraph}. For a demographic $R$ (e.g.~Hispanic), we define a filtration on the graph as $f(v) = 1 - R(v)$ where $R(v) \in [0,1]$ is the proportion of demographic $R$ in tract $v$. Thus for a given threshold $t \in [0,1]$, the sublevel set $X_t$ consists of all tracts with $R$ population share $R(v) \geq 1-t$. The sequence of sublevel sets gives rise to a persistence diagram (Section \ref{sec:background}). At least one death value in each persistence diagram is $\infty$ since when $t = 1$ we have at least one connected component, and we replace these values by $1$ in our diagrams so that we can compute finite Wasserstein distances necessary for our analysis. 

We have now fully described how to obtain a persistence diagram from a city and choice of demographic $R$. We can also give the following non-technical description of the persistence diagram for a region, wherein each point in the persistence diagram has a specific meaning. See Figure \ref{fig:chiPD} for a visual example.
\begin{itemize}
    \item Each point $(b,d)$ in the persistence diagram represents a census tract which has higher $R$ population share than any of its neighboring tracts (the \emph{birth tract}).
    \item The value $b$ is the non-$R$ population share of the birth tract.
    \item We can describe the meaning of the death value $d$ using paths. In particular, $d$ is the largest number which satisfies the following statement: for every path from the birth tract which ends at a tract of higher $R$ share, the path must pass through a tract of non-$R$ share at least $d$.
\end{itemize}

As we will soon see, persistent points can correspond to a wide range of phenomena ranging from individual tracts with artifically high Black or Hispanic populations to large regions of high Black or Hispanic percentage. It will be convenient to refer to all of these by the general term \emph{cluster}. In other words: each high persistence point corresponds to a cluster.

\begin{remark}
We could also employ $1^{st}$-dimensional persistent homology to detect ``holes'' in the data as was done in \cite{feng2021persistent}. Specifically, one could look for holes in the non-Black or non-Hispanic population in order to find encircled clusters of high Black or Hispanic population. However, this does not add anything new as these clusters are already detected by the $0^{th}$-dimensional persistence diagram, except for clusters on the boundary (which we think should not be ignored). Thus, we stick to the more readily explainable $0^{th}$-dimensional analysis in this paper.
\end{remark}

\begin{remark}
The ability to connect points in a persistence diagrams to specific locations which can be further studied is a key feature of our method. In some sense it is similar to how local Moran's I and other LISA indicators~\cite{anselin1995local} go beyond a single number by detecting  specific locations with high spatial auto-correlation. 
\end{remark}

\section{Persistence diagram of a single city: Chicago}\label{sec:single}

\subsection{Persistence diagram}

\begin{figure}
    \centering
    \includegraphics[width=0.9\textwidth]{ 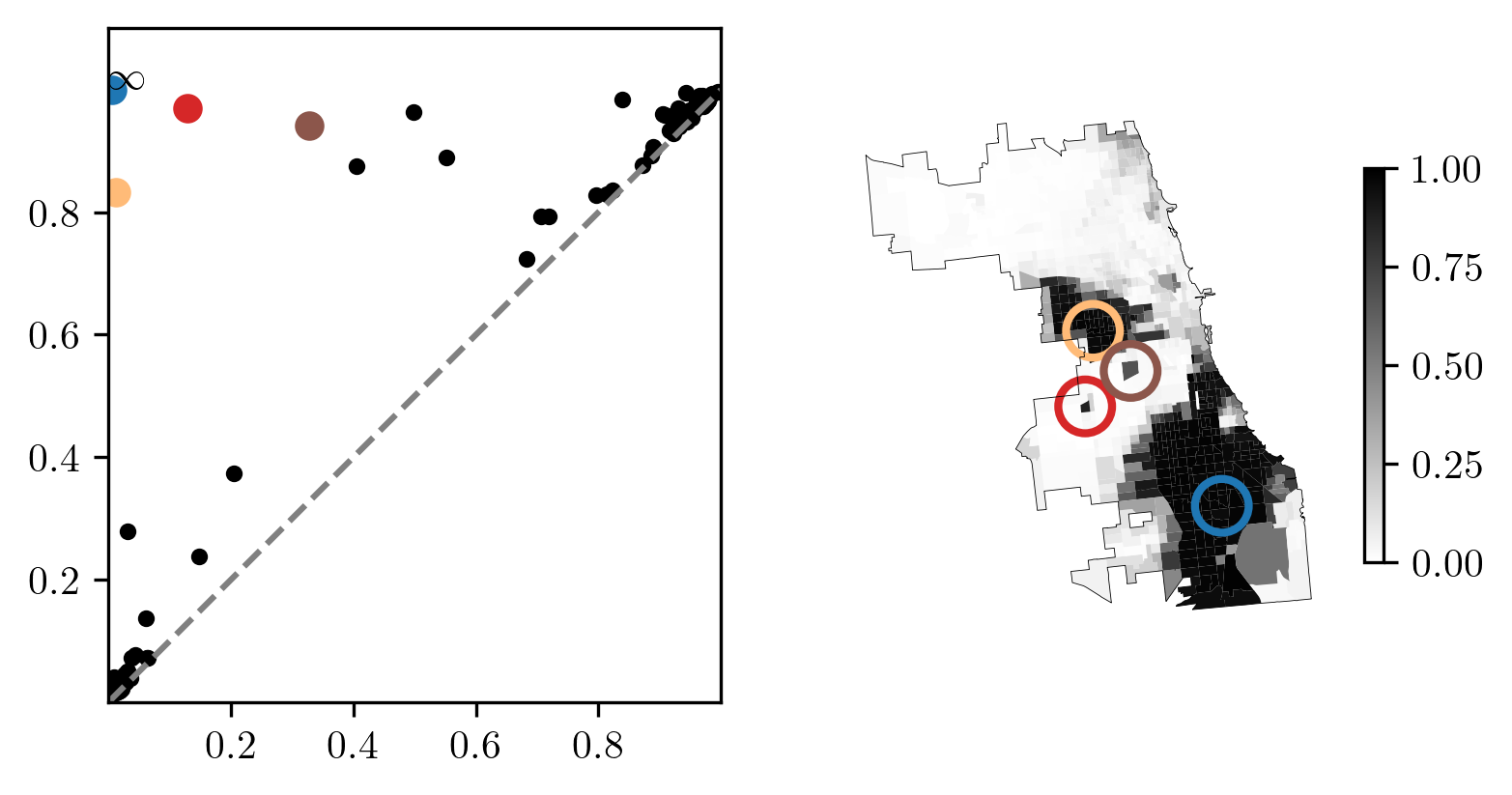}
    \caption{The persistence diagram (left) and choropleth (right) for Black population in Chicago in 2010. The four most persistent points are marked with colors on the persistence diagram and their corresponding birth tracts are circled on the map. These tracts represent large communities on the South and West side, a public housing development, and a correctional facility.}
    \label{fig:chiPD}
\end{figure}

Before we look at an analysis of all 100 cities using various persistence-based visualization techniques, it is worth briefly describing the kind of features that a persistence diagram detects. We use as our example the city of Chicago, its Black population, and the year 2010. Figure \ref{fig:chiPD} shows the persistence diagram side by side with a choropleth showing Black population percentage (as a fraction between 0 and 1). Each point in the persistence diagram corresponds to a tract with higher Black percentage than its neighbors (see Section \ref{sec:constructing}); for the four most persistent points we have circled those tracts in Figure \ref{fig:chiPD} with the corresponding color. Notice how these tracts are separated from one another by areas of low Black percentage (which is what creates high persistence). The two points near the $y$-axis in blue and beige represent large Black neighborhoods in the South and West sides of the city. The next most persistent point (shown in red) corresponds to just one tract with a census population (in 2010) of 2319 people, 2291 of whom were Black. This tract contains the public housing project of LeClaire Courts. The fourth most persistent point, in brown, again represents a single tract. This single tract contains the Cook County Department of Corrections, one of the largest single-site jails in the United States. According to the Census, the population of this tract was 11309 people, 7603 of them Black, a much higher percentage than nearby tracts. 

\begin{remark}
    This example illustrates that persistent homology does not distinguish between ``natural'' demographic patterns and artificial concentration and segregation arising from the connection between race and factors such as income and incarceration rates. In U.S.~cities, the effects of public policy decisions on demographic distributions is well-studied \cite{massey1988dimensions, wilson2012truly}, so we should always be careful to investigate the reasons behind the particular patterns being detected rather than treating them as natural or inherent properties of a city or region.
\end{remark}

\begin{figure}
    \centering
    \includegraphics[width=0.8\textwidth]{ 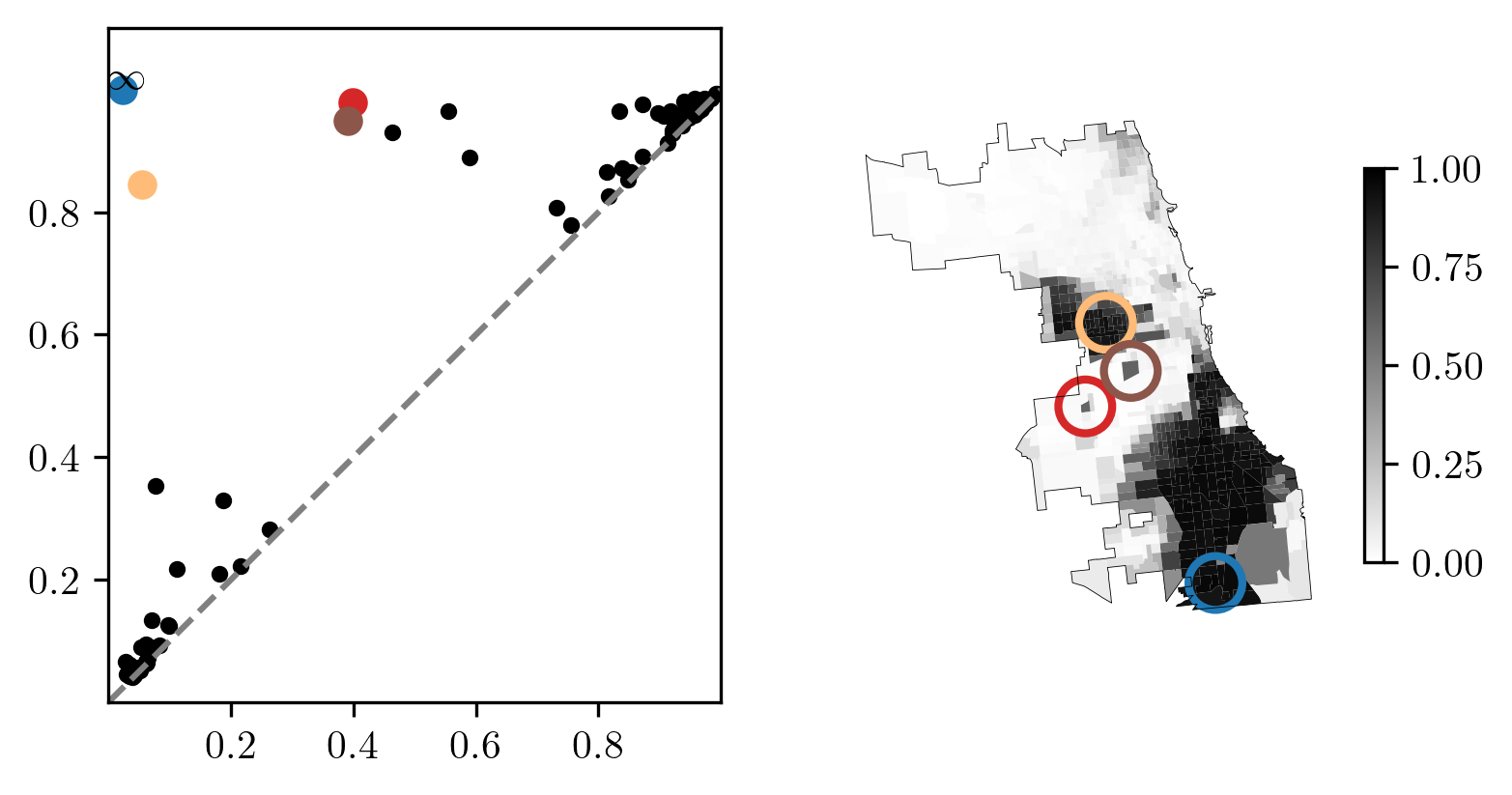}
    \caption{The persistence diagram (left) and choropleth (right) for Black population in Chicago in 2020, to be compared with the 2010 diagram in Figure \ref{fig:chiPD}. }
    \label{fig:chiPDtime}
\end{figure}

\subsection{Distances between persistence diagrams}

We also use Chicago to demonstrate the kind of changes that can be detected by comparing persistence diagrams. For this purpose, we use 2010 and 2020 census data for Chicago to create two persistence diagrams, shown in Figure~\ref{fig:chiPD} and Figure~\ref{fig:chiPDtime} respectively. We see that only one very persistent point has substantially moved, namely the red point. The others have stayed the same, indicating that little has changed about the overall distribution of the Black population of Chicago. Indeed, we will see later on that, as measured by Wasserstein distance, Chicago has changed noticeably less than many cities of similar size. The point that moved represents the tract containing LeClaire Courts public housing project, whose tract went from a Black population of 2291 (out of 2319 total) in 2010 to a Black population of 1029 (out of 1710 total) in 2020. Though restricted to a single tract, this is a significant demographic shift.

We therefore conclude that highly persistent points can represent a variety of different demographic features including large communities or neighborhoods, housing developments and correctional facilities. These features together comprise the way that a particular population is distributed in space. Using persistence we can compare these distributions both between cities and across time.

\section{Finding trends and outliers} \label{sec:scatter}

\subsection{Visualizing differences using dimension reduction}\label{sec:W1outliers}

\begin{figure}
    \centering
    \begin{subfigure}[t]{0.4\textwidth}
        \centering
    Black (2020)
    \includegraphics[width=\textwidth]{ 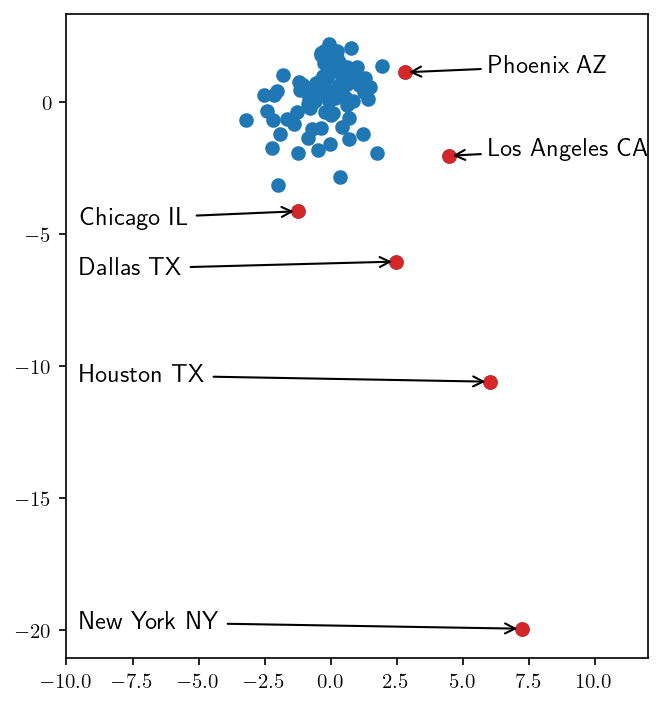}
    \end{subfigure}
    \begin{subfigure}[t]{0.4\textwidth}
    \centering
    Hispanic (2020)
    
    \includegraphics[width=\textwidth]{ 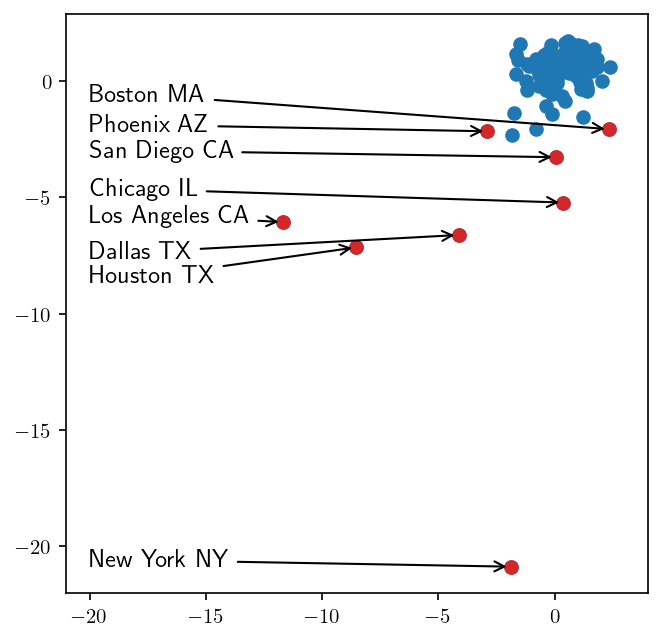}
    \end{subfigure}
    \caption{2D visualizations of the relationships between cities as measured by persistent homology created using multi-dimensional scaling. Red points indicate outliers.}
    \label{fig:mds}
\end{figure}

Let $C_1,\ldots, C_{100}$ denote the 100 cities under consideration. Suppose we have fixed a choice of demographic group (Black or Hispanic) and a year (2010 or 2020). For any pair of cities $C_i$, $C_j$, we can compute a Wasserstein distance between them by converting them each into persistence diagrams (Section \ref{sec:constructing}) and computing the $1$-Wasserstein distance between these diagrams. What we obtain is a $100 \times 100$ distance matrix. We can then use a dimensionality-reduction method to embed the 100 cities as 100 points in the plane while maintaining these distances as best as possible. We use multi-dimensional scaling (MDS) to perform the embedding, with the results shown in Figure~\ref{fig:mds}. 

The basic structure observed in these plots is of a central cluster surrounded by a few outlier points. We use the Local Outlier Factor (LOF) algorithm~\cite{breunig2000lof} on the original distance matrix to find all outliers, which we show as red points and list in Table \ref{tab:outliers}. 

\begin{remark}
Computing the Local Outlier Factor (LOF) requires the choice of a ``number of neighbors'' parameter $k$ and a threshold $\varepsilon \geq 1$. There is no rigorous way to choose these values, because there is no formal definition of ``outlier''. In our analysis, we compute LOF for the range $k = 10, 11 \ldots, 19$, average the results, and pick a threshold of $\varepsilon = 2.$
\end{remark}

\begin{table}[h]
    \centering
    \begin{tabular}{|l|l|l|}
    \hline
        Black (2020) & Hispanic (2020) & pop.~rank \\
        \hline
        New York NY & New York NY & 1 \\
        Los Angeles CA & Los Angeles CA & 2 \\
        Chicago IL & Chicago IL & 3 \\
        Houston TX & Houston TX & 4 \\
        Phoenix AZ & Phoenix AZ & 6 \\
         & San Diego CA & 8 \\
        Dallas TX & Dallas TX & 9 \\
        & Boston MA & 23 \\
        \hline
    \end{tabular}
    \caption{Outliers based on Wasserstein distance listed for Black (2020) and Hispanic (2020) data. The ranks (based on total population) are shown, indicating that the largest cities are more likely to be outliers.}
    \label{tab:outliers}
\end{table}

We note that outliers tend to be the larger cities. Note that persistence diagrams do not depend on city size since they do not deal with raw population counts, only percentages. Indeed, a city which is divided into a 100\% Black side and a 0\% Black side would give the same persistence diagram whether it was large or small. Thus, the preponderance of very large cities among outliers in our study indicates that the distribution of Black and Hispanic people in a very large city is not merely a scaled-up version of their distribution in a typically smaller city. We further investigate what makes a city an outlier in the next section. 

\subsection{Total persistence}

\begin{figure}[h]
    \centering
    \begin{subfigure}{0.4\textwidth}
    \centering
    Black (2020)
    \includegraphics[width=\textwidth]{ 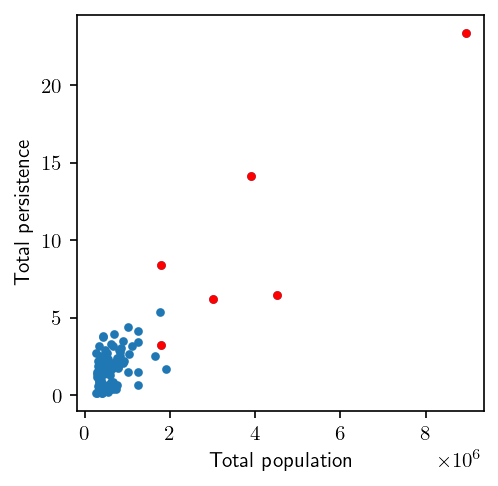}
    \end{subfigure}
    \begin{subfigure}{0.4\textwidth}
    \centering
    Hispanic (2020)
    \includegraphics[width=\textwidth]{ 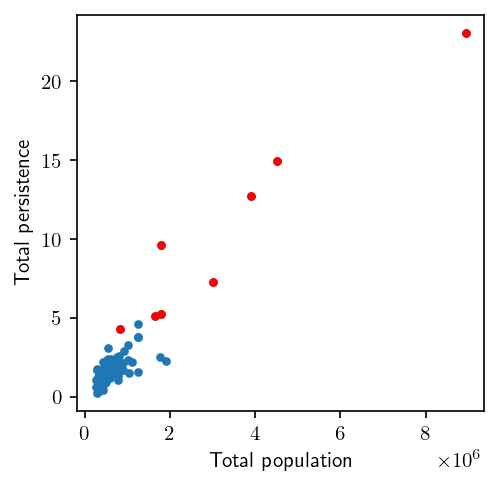}
    \end{subfigure}
    
    \begin{subfigure}{0.4\textwidth}
    \centering
    \includegraphics[width=\textwidth]{ 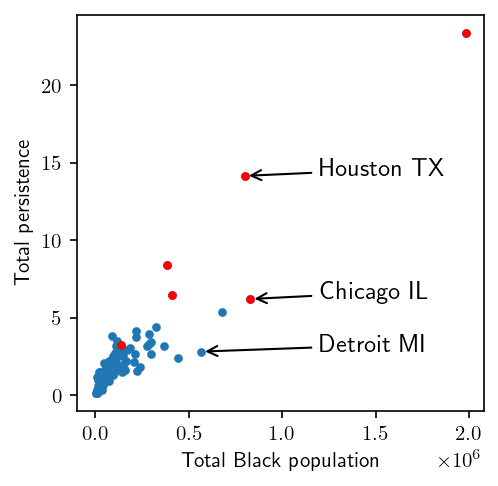}
    \end{subfigure}
    \begin{subfigure}{0.4\textwidth}
    \centering
    \includegraphics[width=\textwidth]{ 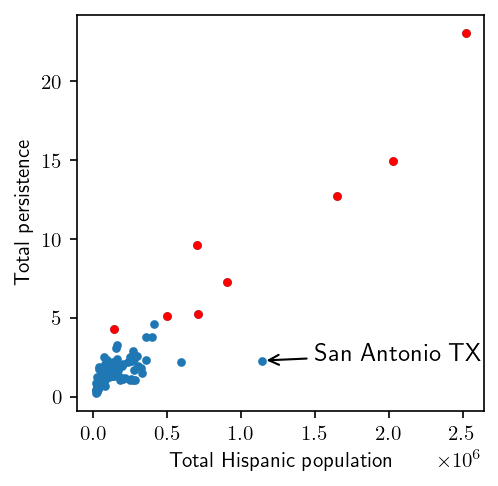}
    \end{subfigure}
    
    \caption{Scatter plots of total population (above) and total Black/Hispanic population (below) vs total persistence. Red points indicate outliers. }
    \label{fig:tps}
\end{figure}

To further examine what makes a city an outlier, we compute the total persistence (Definition \ref{def:tp}) of each city and compare it to the total population and the Black/Hispanic population -- see Figure \ref{fig:tps}. The correlation between population size and total persistence is clear in all cases ($r > 0.9$). 

Most of the outliers tend to be cities with high populations and high total persistence. High total persistence typically arises from having many distinct neighborhoods with high Black/Hispanic population. This kind of structure therefore appears to be more common in large cities where there is enough space for these neighborhoods to develop. We briefly take a closer look at some specific cities labelled in Figure \ref{fig:tps}.

\subsubsection{Detroit MI and San Antonio TX} 
\begin{figure}
    \centering
    \begin{subfigure}{0.48\textwidth}
    \centering
    San Antonio (Hispanic)
          \includegraphics[width=\textwidth]{ 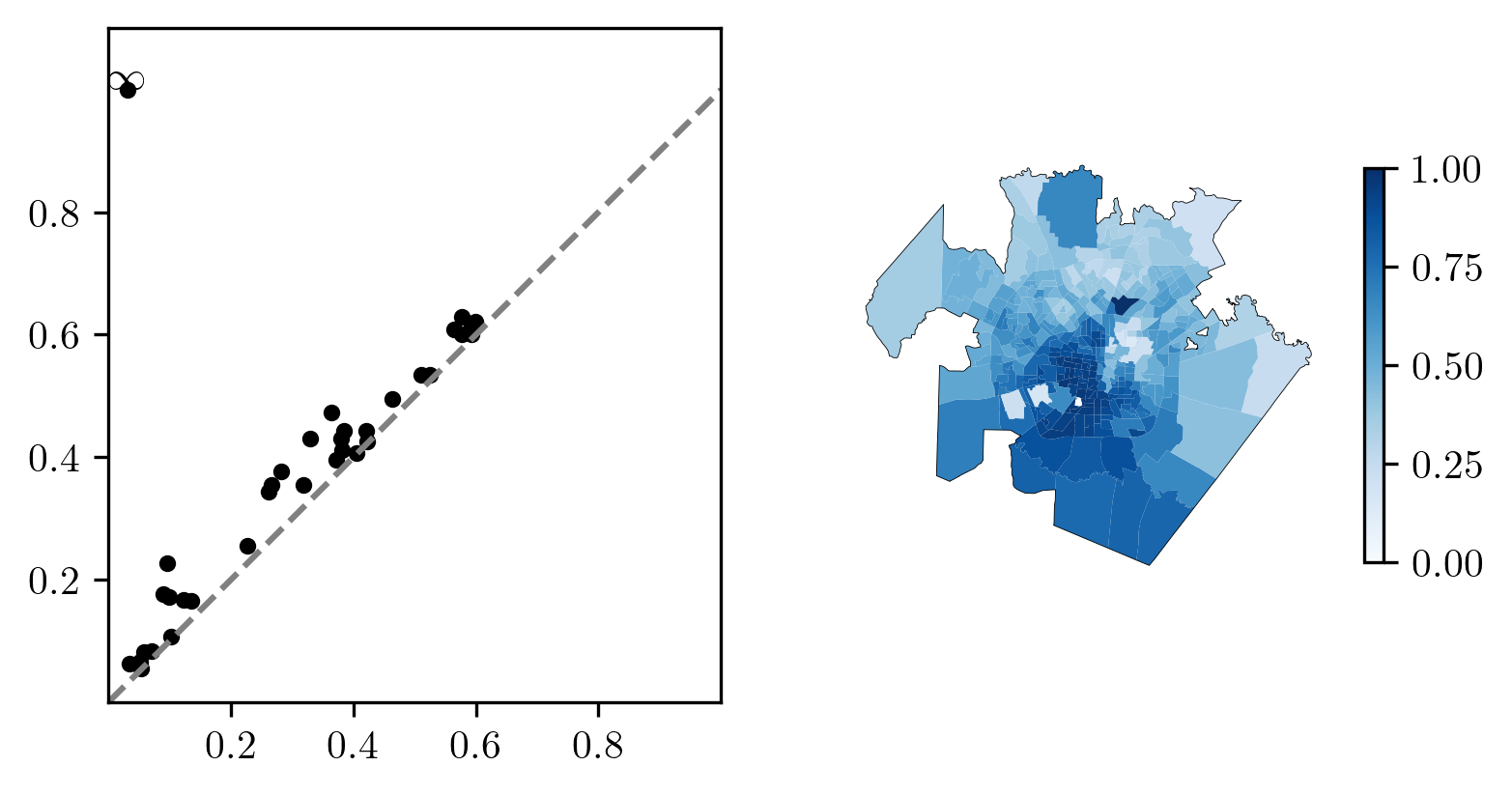}
    \end{subfigure}
    \begin{subfigure}{0.48\textwidth}
    \centering
   Detroit (Black)
          \includegraphics[width=\textwidth]{ 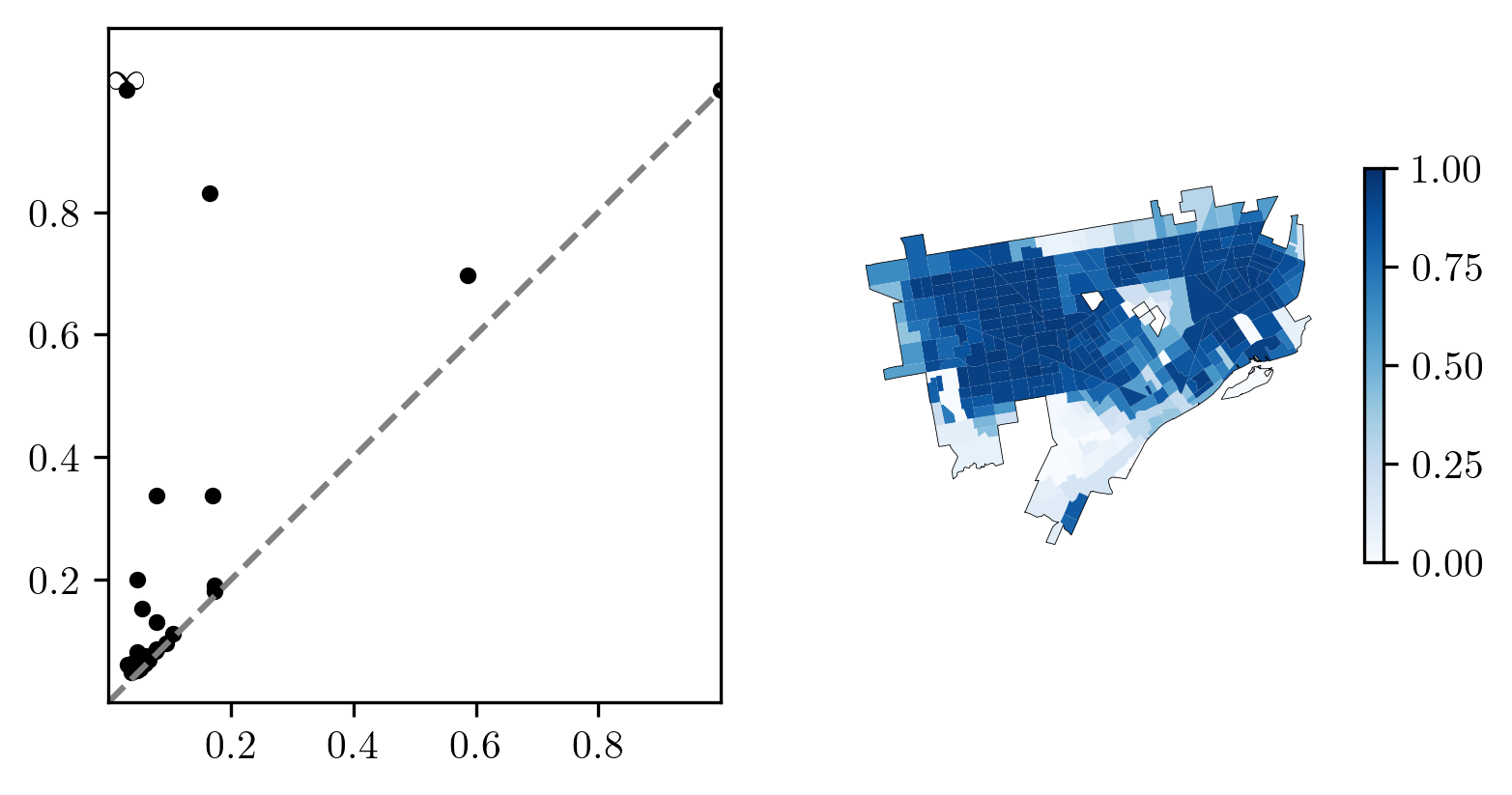}
    \end{subfigure}
    \caption{San Antonio TX showing Hispanic \% and Detroit MI showing Black \%, with the corresponding persistence diagrams.}
    \label{fig:sanantonio}
\end{figure}

In Figure \ref{fig:tps} we can see that Detroit MI and San Antonio TX have lower total persistence than one might expect based on other cities of similar population sizes -- for Black population in the case of Detroit and for Hispanic population in the case of San Antonio. In Figure \ref{fig:sanantonio}, we see that each of these cities has very few high persistence points (in fact San Antonio only has one). This means that almost all tracts with high Black or Hispanic population are restricted to a single large cluster. This could represent evidence of high levels of segregation in these cities, compared to cities of similar size. These large cities are also not outliers, which tells us they may have more in common with smaller cities than with the largest cities in the dataset which have a complicated neighborhood structure.

\subsubsection{Chicago vs Houston} Despite having similar Black population numbers, Chicago and Houston differ greatly in their total persistence for Black percentage. We contrast the two persistence diagrams in Figure \ref{fig:chicagohouston} below. We can clearly see that Houston has many more points away from the diagonal than Chicago, which indicates that Houston's predominantly Black neighborhoods are more numerous and spread out. By contrast, Chicago in 2020 had two large Black neighborhoods, and almost no neighborhoods which were 60-80\% Black at their core.

\subsubsection{Phoenix AZ} 
\begin{figure}
    \centering
    \begin{subfigure}{0.6\textwidth}
            \includegraphics[width=\textwidth]{ 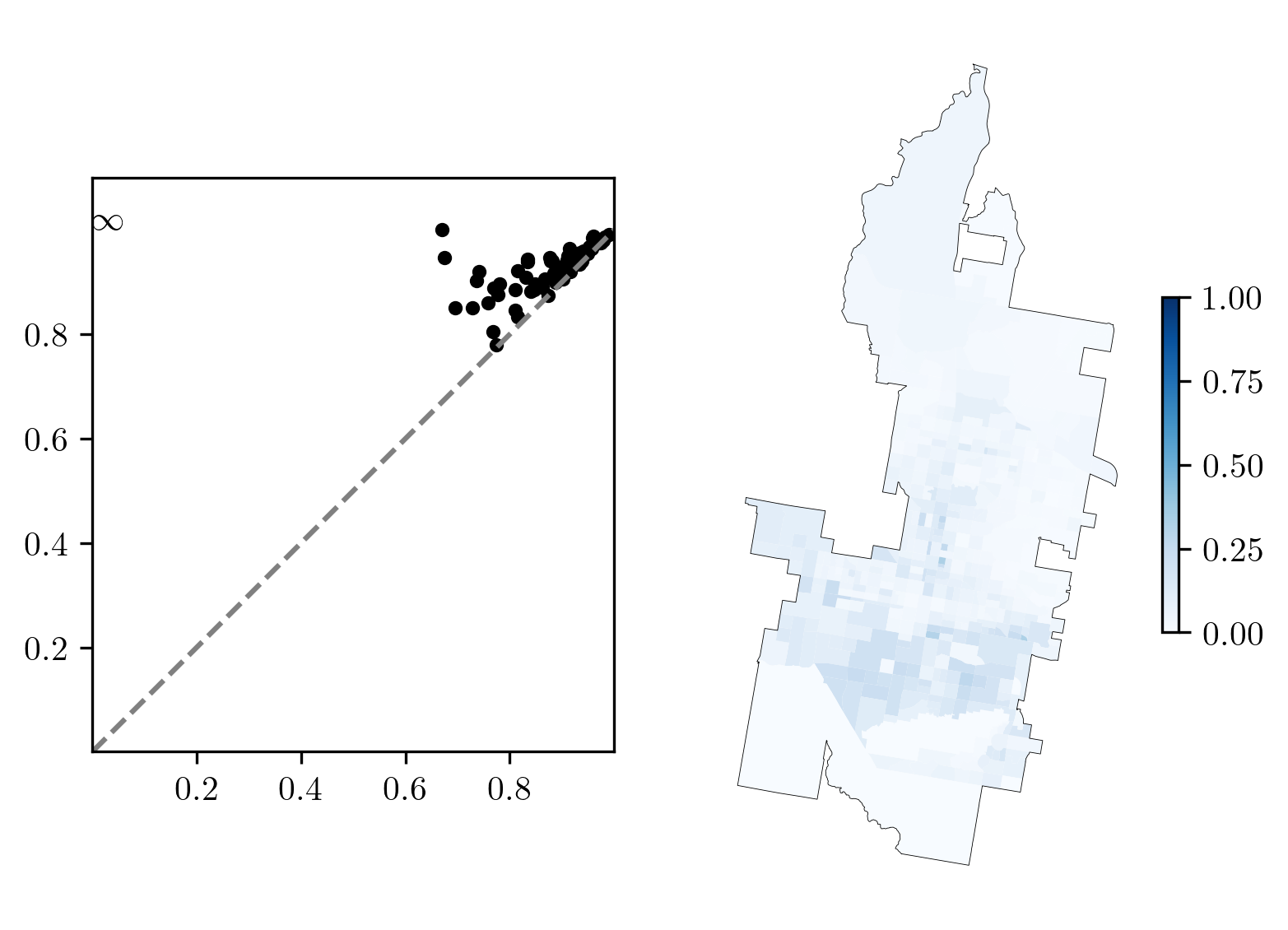}
    \end{subfigure}
    \caption{Phoenix AZ showing Black \% and persistence diagram.}
    \label{fig:phoenixaz}
\end{figure}

We close with an example of an outlier which is not easily explained by total persistence alone: Phoenix AZ for Black percentage. Indeed, Phoenix's population is only about 7\% Black, and the highest Black percentage of any tract is 33\%, very different from either Chicago or Houston. What makes Phoenix unique is the large number of points representing clusters anchored by tracts which are approximately 25\% Black. These are the points away from the diagonal in the upper right corner of the persistence diagram in Figure~\ref{fig:phoenixaz}. In other words, cities with the same low percentage of Black population tend not to have as many distinct mixed neighborhoods as Phoenix.

\begin{figure}
    \centering
    \begin{subfigure}{0.45\textwidth}
    \centering
    Chicago
            \includegraphics[width=\textwidth]{ 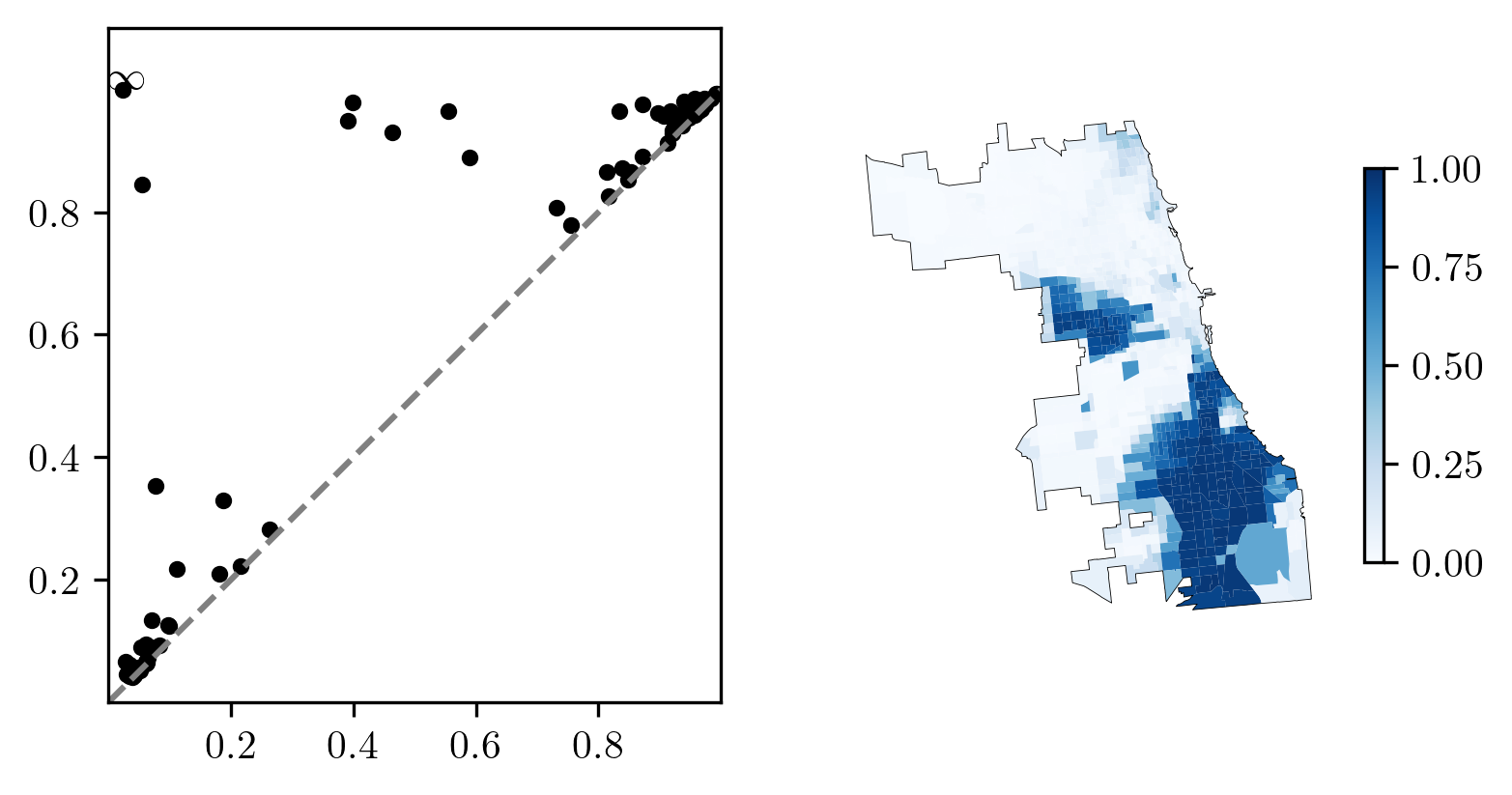}
    \end{subfigure}\hfill 
    \begin{subfigure}{0.45\textwidth}
    \centering
    Houston
            \includegraphics[width=\textwidth]{ 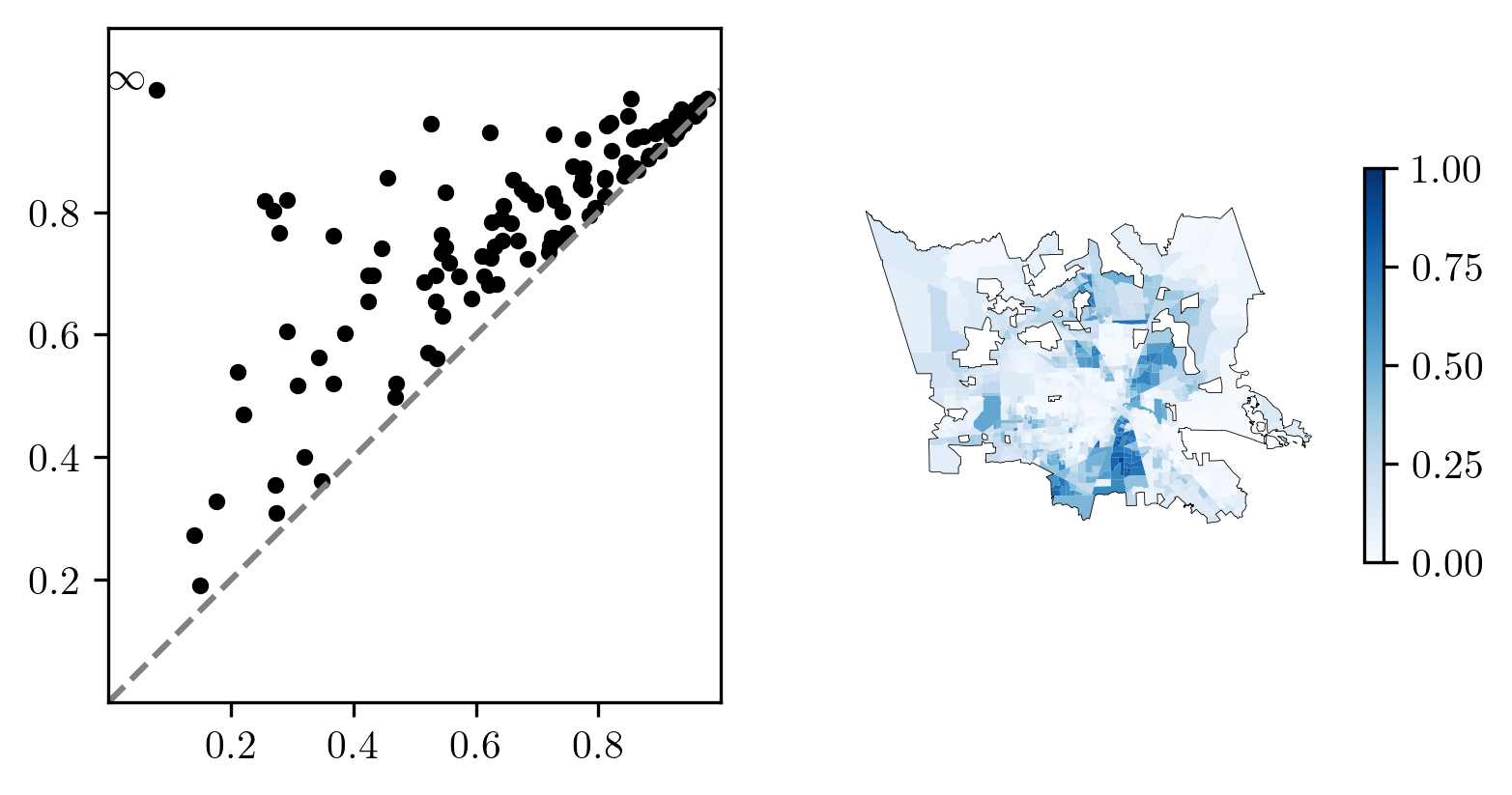}
    \end{subfigure}
    \caption{Contrasting Houston and Chicago with respect to  Black population using choropleths and persistence diagrams.}
    \label{fig:chicagohouston}
\end{figure}

We also examine whether total persistence is related to classical measures of segregation, choosing for this purpose one of the most common measures of segregation, the \emph{dissimilarity index} (see \cite{massey1988dimensions} for a definition and more examples of segregation measures). We compute Black vs non-Black and Hispanic vs non-Hispanic dissimilarity indices for all cities in our dataset.

\begin{figure}
    \centering
    \begin{subfigure}{0.36\textwidth}
    \centering
    Black (2020)
        \includegraphics[width=\textwidth]{ 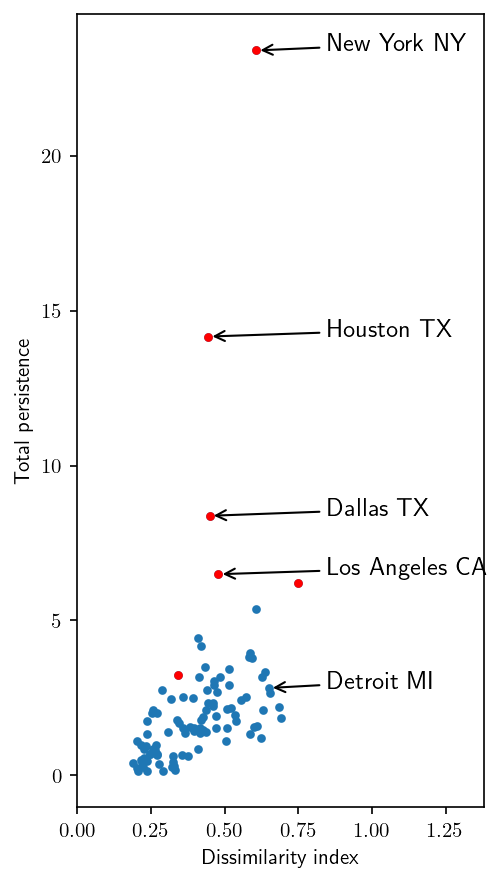}
    \end{subfigure} \hspace{1cm}
    \begin{subfigure}{0.36\textwidth}
    \centering
    Hispanic (2020)
        \includegraphics[width=\textwidth]{ 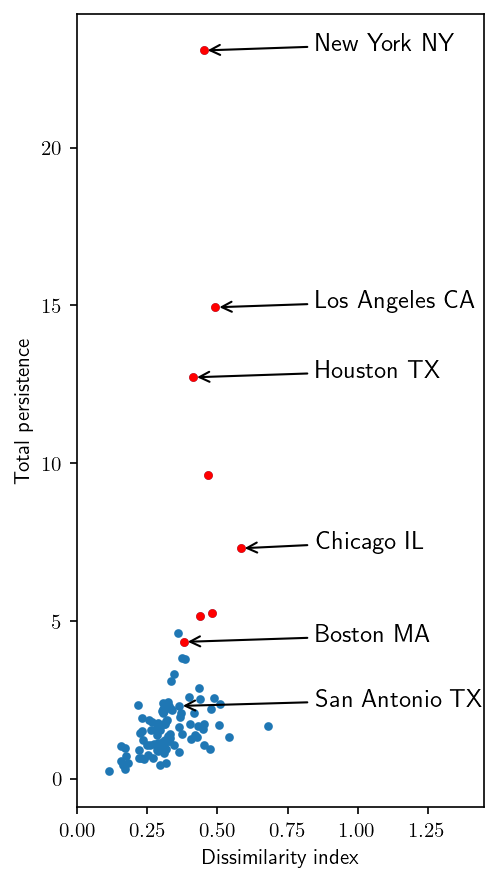}
    \end{subfigure}
    \caption{Plotting dissimilarity index (a measure of segregation) against total persistence. We see a positive but weak correlation. Outliers from Figure \ref{fig:mds} are shown in red.}
    \label{fig:dis}
\end{figure}

Figure \ref{fig:dis} shows dissimilarity index (DI) against total persistence. For both Black and Hispanic, we find a correlation coefficient of about $0.39$, indicating a positive but weak correlation. Indeed, we see that while there is a slight upward trend, for each DI value there is a range of possible total persistence values. 

Note that DI does not take the spatial structure of tracts into account; it operates purely on list of population totals by tract. Thus while DI can detect something about the demographic distributions by looking at a city tract-by-tract, total persistence brings in another dimension. In particular, DI does not distinguish between high Black/Hispanic percentage tracts being clustered into one place or being spread apart across the city; while our persistence diagrams are specifically designed to detect this distinction. We therefore propose that the combination of a measure like DI with persistence measures is a good indicator of the demographic structure of a city, and can tell us a lot about the nature and scale of segregation. 

\section{Measuring change over time}\label{sec:change}

For each city $C_i$ we consider the persistence diagrams generated from 2010 and 2020 data respectively, and compute the $1$-Wasserstein distance between them. Recall that the Wasserstein distance measure is an additive measure, so that this distance represents the total change to all persistence points under the best possible matching. The ability to even compute this kind of distance is an advantage of our persistence-based approach. Indeed, census tracts can and do change from one cycle to the next, so it is impossible to directly track percentage changes at the tract level. Our method does not require us to track individual tracts, and so avoids this issue.

\begin{figure}
    \centering
    \begin{subfigure}{0.4\textwidth}
        \includegraphics[width=\textwidth]{ 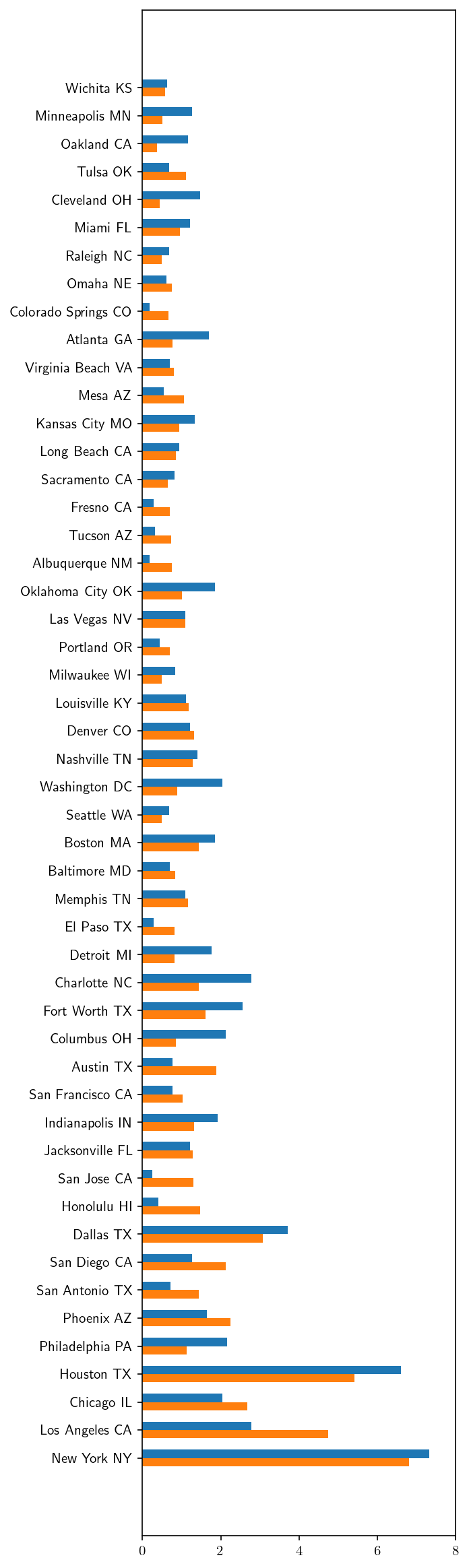}
    \end{subfigure}%
     \begin{subfigure}{0.4\textwidth}
        \includegraphics[width=\textwidth]{ 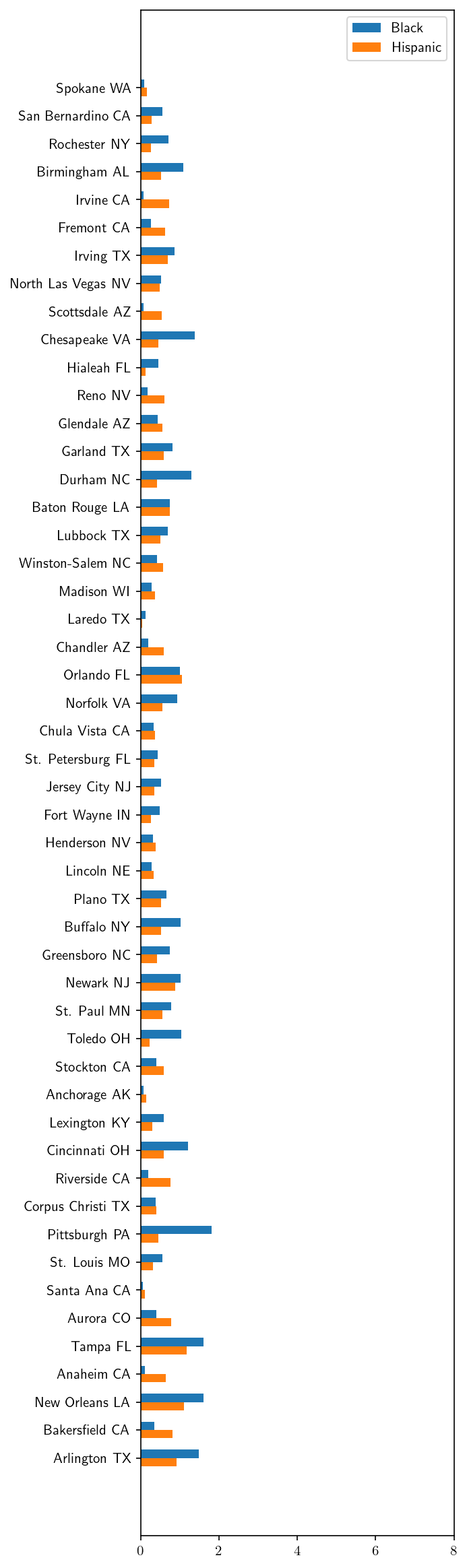}
    \end{subfigure}
    \caption{Wasserstein distances between 2010 and 2020 data for each city.}
    \label{fig:wassovertime}
\end{figure}

\begin{figure}
    \centering
    \begin{subfigure}{0.45\textwidth}
    \centering
    Pittsburgh (2010)
            \includegraphics[width=\textwidth]{ 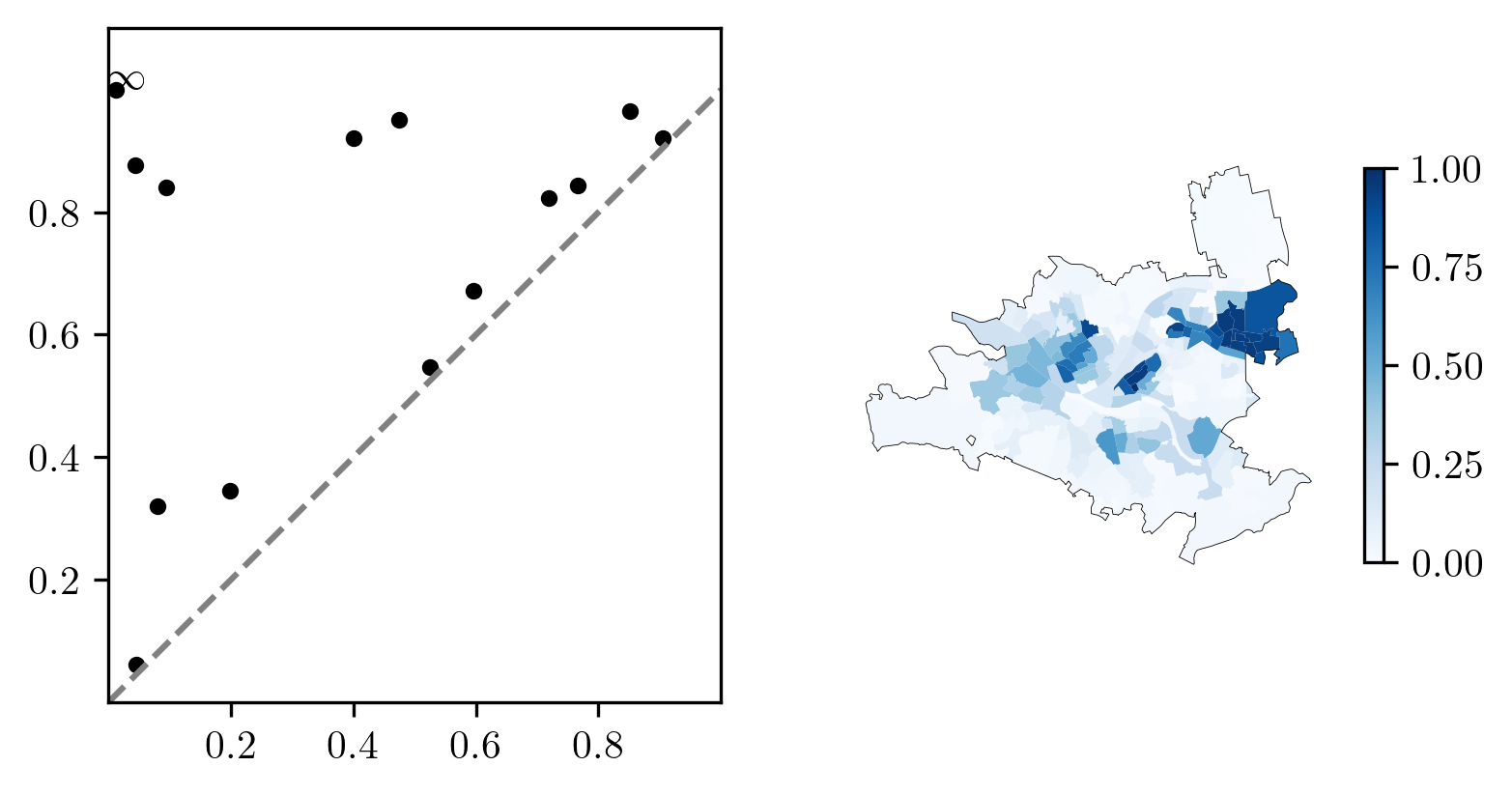}
    \end{subfigure}\hfill 
    \begin{subfigure}{0.45\textwidth}
    \centering
    Pittsburgh (2020)
            \includegraphics[width=\textwidth]{ 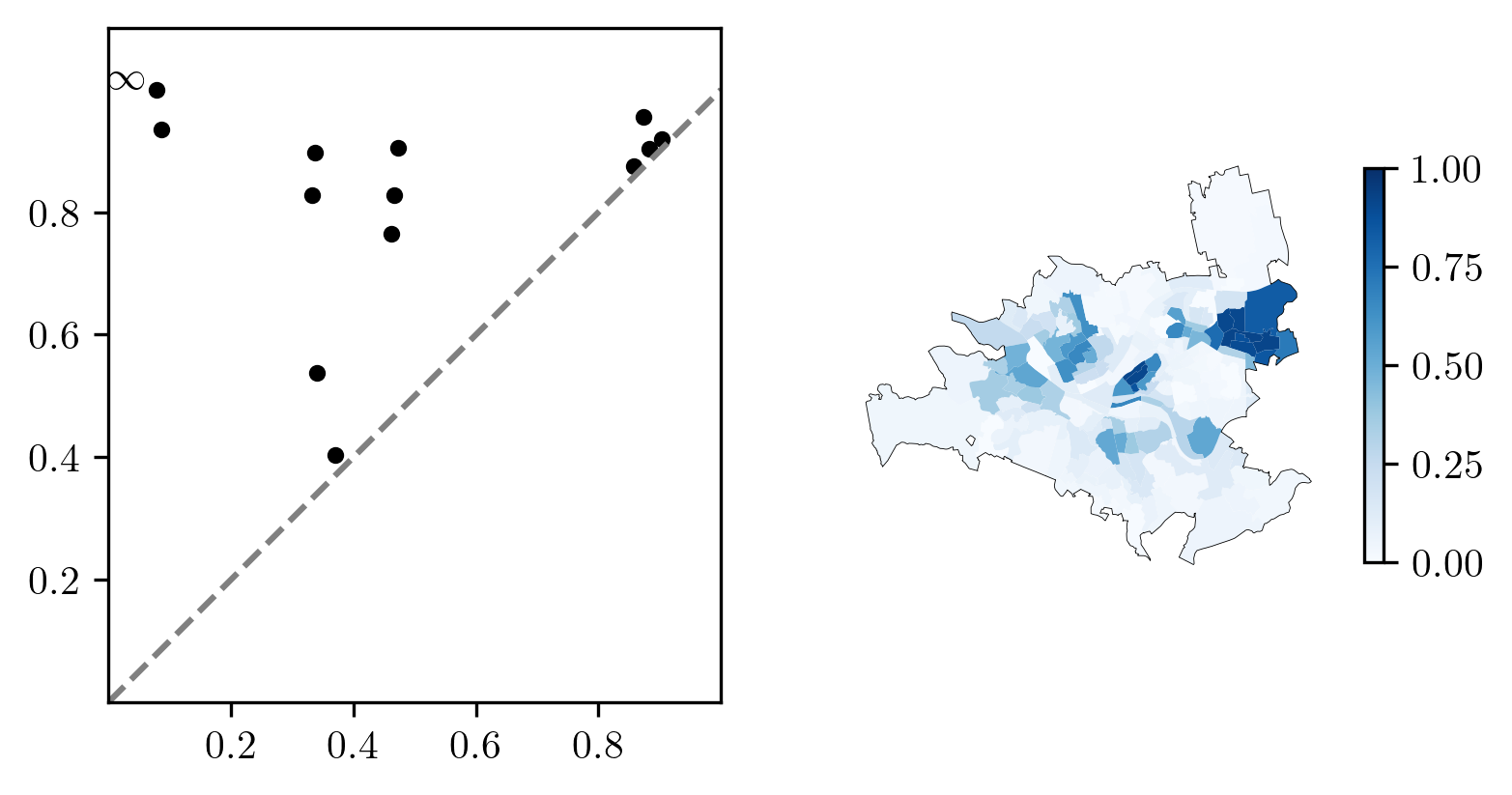}
    \end{subfigure}
    \caption{Contrasting Pittsburgh PA in 2010 and in 2020 using Black population.}
    \label{fig:pittsburgh}
\end{figure}

The Wasserstein distance between the 2010 and 2020 persistence diagrams for each city are shown in Figure \ref{fig:wassovertime}. We see that the largest cities typically exhibit the greatest change as measured by Wasserstein distance. This is expected, as more neighborhoods allows for more total change. Some examples of cities which stand out are Pittsburgh PA which has a larger change for Black population than most cities of similar size. This change is clearly visible on the maps and persistence diagrams in Figure \ref{fig:pittsburgh}. We note that Chicago has a low amount of change compared to other cities of similar or smaller size sich as Houston TX or Dallas, suggesting that the fundamental structure of Chicago's demographics has not changed much in ten years.

\section{Clustering analysis} \label{sec:clustering}

\begin{figure}[h]
    \centering
    \begin{subfigure}{0.48\textwidth}
        \includegraphics[width=\textwidth]{ 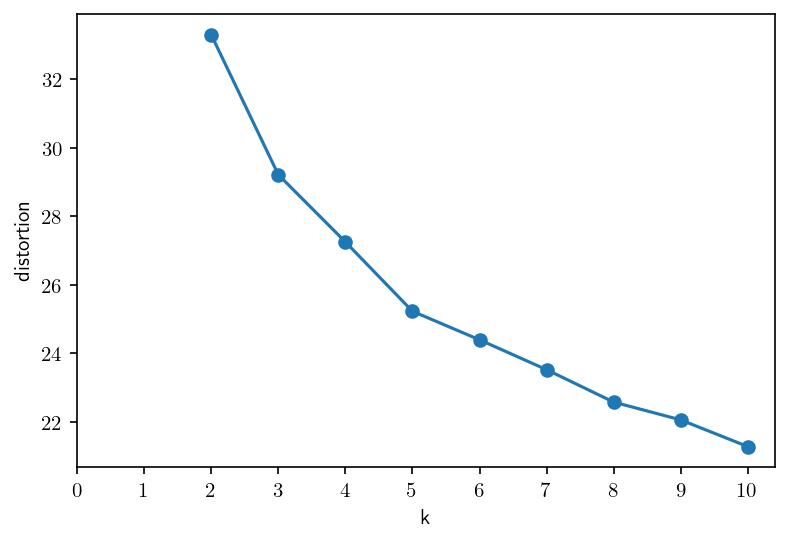}
        \caption{Black \%}
    \end{subfigure}
    \begin{subfigure}{0.48\textwidth}
        \includegraphics[width=\textwidth]{ 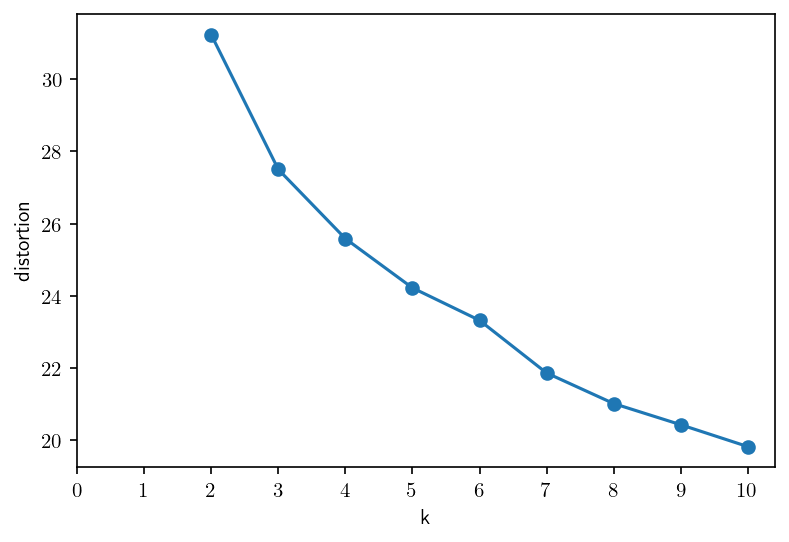}
        \caption{Hispanic \%}
    \end{subfigure}
    \caption{Elbow plots for the $k$-means clustering of persistence diagrams.}
    \label{fig:elbow}
\end{figure}

In order to gain a better picture of how the 100 cities compare to one another, we cluster the persistence diagrams for the 100 cities (with 2020 data) using the $k$-means clustering algorithm introduced in \cite{marchese2017k}. This algorithm takes the number of clusters $k$ as a parameter and returns $k$ clusters as well as a Fr\'{e}chet mean (which we will call the \emph{cluster mean}) for each cluster. A Fr\'{e}chet mean is a locally optimal ``average'' diagram initially introduced in \cite{turner2014frechet}. The cluster labels give us a classification of cities by their persistence data, while the cluster means allow us to characterize the common characteristics within each cluster. To choose the value of $k$, we adapt a standard procedure from classical $k$-means called an \emph{elbow plot}, shown in Figure~\ref{fig:elbow}. For a range of values of $k$, we compute the distortion
$$
D = \sum_i W_2(P_i, M_{\sigma(i)})^2
$$
where $P_i$ is the $i^{th}$ city's persistence diagram, $M_j$ is the $j^{th}$ cluster mean, and $\sigma(i)$ is the cluster label assigned to city $i$. We choose a value of $k$ where the distortion appears to begin a roughly linear descent with respect to $k$. The elbow plots in Figure~\ref{fig:elbow} do not present a strong signal for number of clusters, but we chose $k=5$ for both as a point where the graphs start to appear more linear.

For each $k$, we also require an initial set of $k$ cities to act as cluster centers in the first iteration of the algorithm. In order to choose these initial centers effectively, we adapt the $k$-means++ algorithm to the persistence diagram setting \cite{arthur2007k}.

\subsection*{Black population.} The five clusters are shown in Figure \ref{fig:BLACK_map}. We see one cluster with just two cities, New York NY and Houston TX, which have a large number of high persistence points. These two cities had the highest total persistence by some margin (Figure \ref{fig:tps}) and so are naturally grouped into a cluster of their own. Cluster 2 contains cities with a few clusters of high Black population (e.g. Chicago IL, Pittsburgh PA), while cluster 3 contains cities with only one cluster (e.g. Milwaukee WI, Detroit MI). Clusters 4 and 5 consist of cities with lower or less concentrated Black populations.

\subsection*{Hispanic population.} 
\begin{figure}[h]
    \centering
    \begin{tikzpicture}[scale=1.12]
    \node at (0,0) {\includegraphics[width=0.8\textwidth]{ 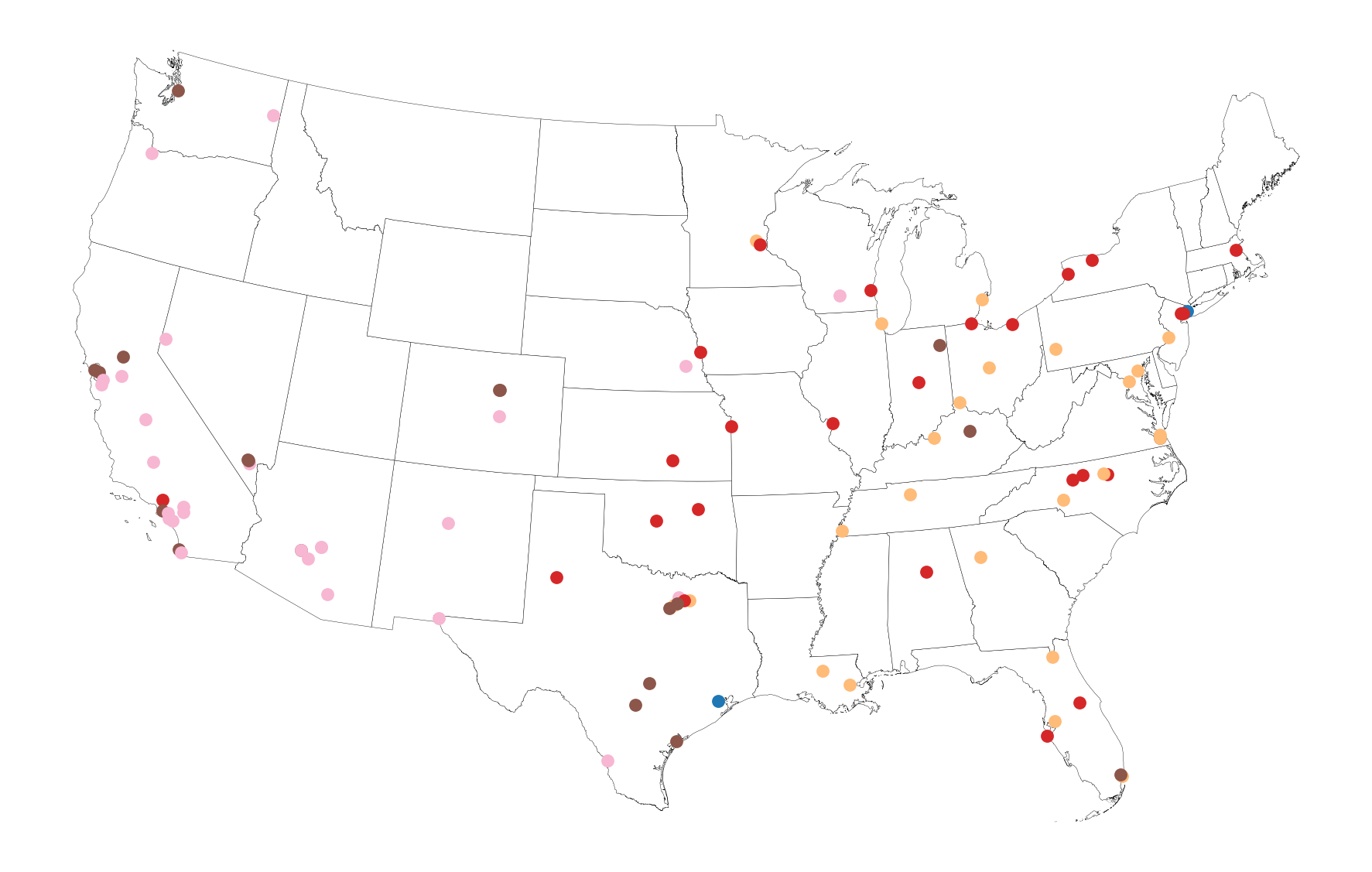}};
    \node at (-5, -4.7) {\includegraphics[width=70pt]{ 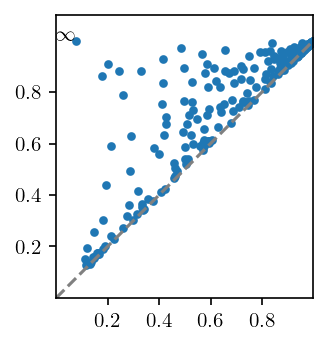}};
    \node at (-3, -4.7) {\includegraphics[width=70pt]{ 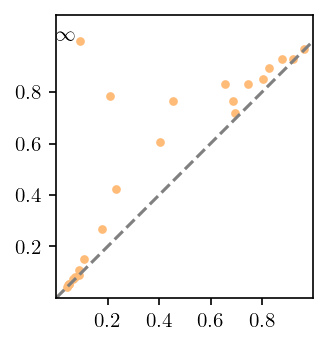}};
    \node at (-1, -4.7) {\includegraphics[width=70pt]{ 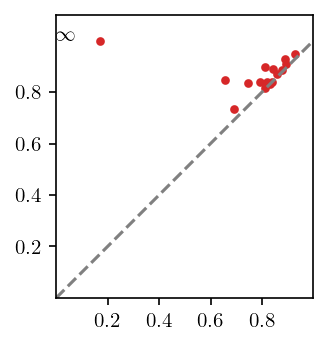}};
    \node at (1, -4.7) {\includegraphics[width=70pt]{ 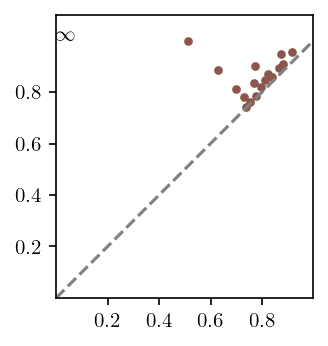}};
    \node at (3, -4.7) {\includegraphics[width=70pt]{ 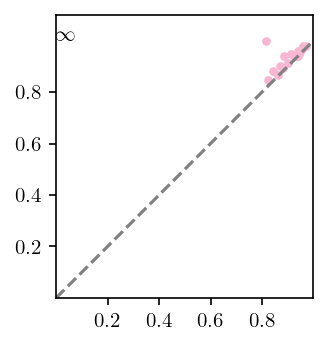}};
    \node at (0, -7.3) {
    \begin{tabular}{c|c|r|r|r|r}
cluster & \#cities & ave. pop & ave. \# tracts & ave. \% BLACK \\ 
1 \textcolor {mycolor1}{$\blacksquare$} & 2 & 6415262.50  & 1652.00  & 21.34 \% \\ 
2 \textcolor {mycolor2}{$\blacksquare$} & 24 & 857486.29  & 237.38  & 33.85 \% \\ 
3 \textcolor {mycolor3}{$\blacksquare$} & 25 & 668597.60  & 190.56  & 24.38 \% \\ 
4 \textcolor {mycolor4}{$\blacksquare$} & 20 & 772336.50  & 181.00  & 11.65 \% \\ 
5 \textcolor {mycolor5}{$\blacksquare$} & 29 & 502276.62  & 116.76  & 4.84 \% \\ 
    \end{tabular}
    };
    \end{tikzpicture}
    \caption{One hundred U.S. cities clustered into $k = 5$ clusters based on their Black \% persistence diagrams. Above, city locations are plotted above and colored by cluster assignment. In the middle, the cluster mean for each cluster is shown. Each cluster mean is an average of the persistence diagrams for that cluster. At the bottom is a table collecting various statistics for the cities in each cluster.}
    \label{fig:BLACK_map}
\end{figure}

\begin{figure}[h]
    \centering
    \begin{tikzpicture}[scale=1.15]
    \node at (0,0) {\includegraphics[width=0.8\textwidth]{ 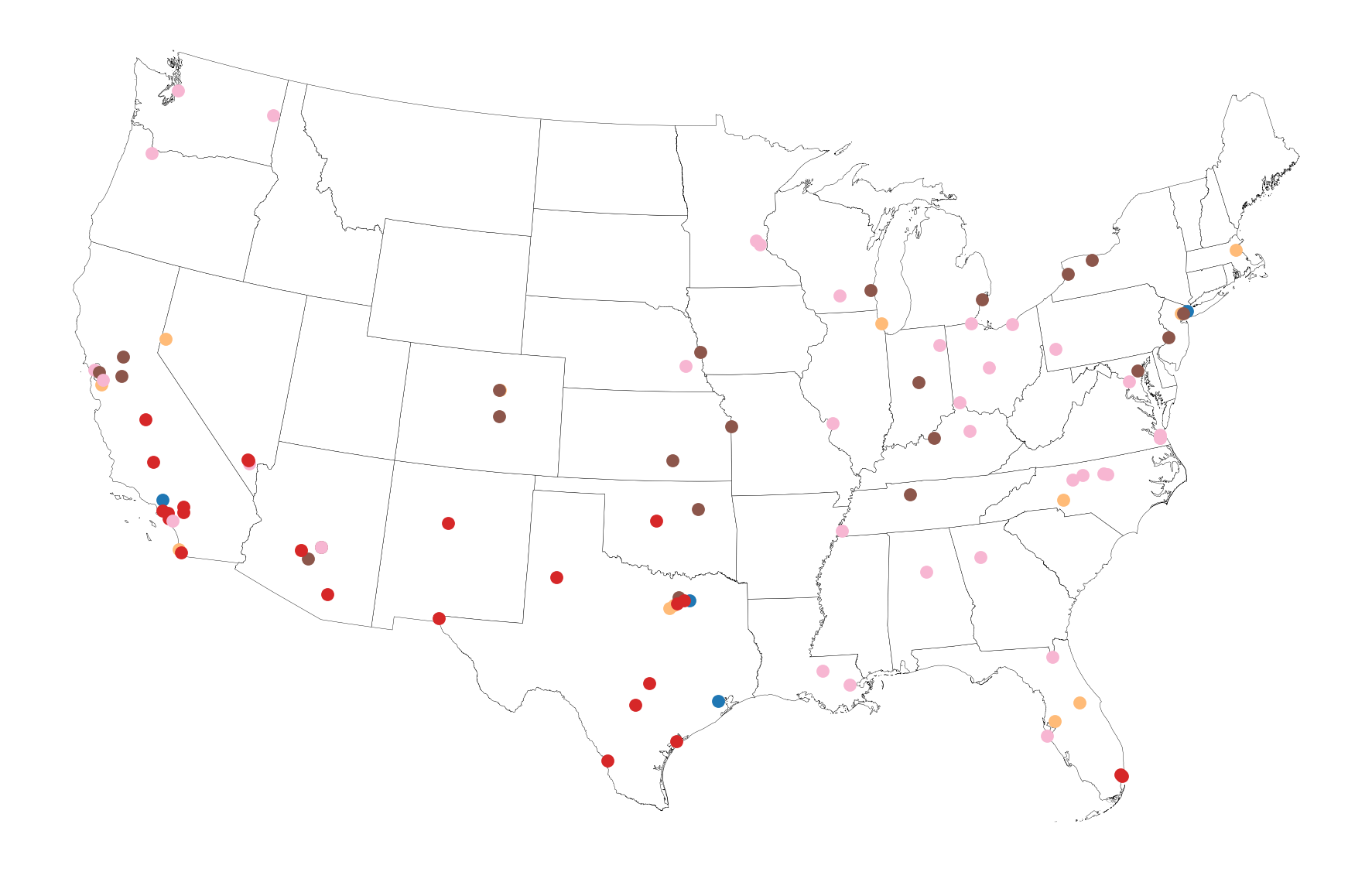}};
    \node at (-5, -4.7) {\includegraphics[width=70pt]{ 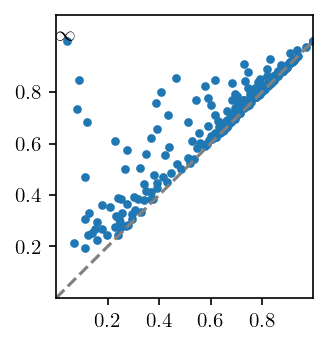}};
    \node at (-3, -4.7) {\includegraphics[width=70pt]{ 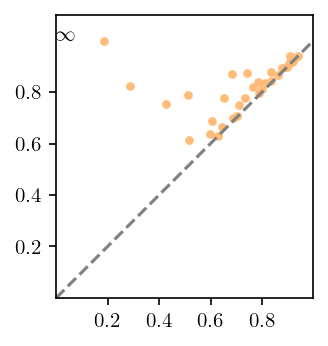}};
    \node at (-1, -4.7) {\includegraphics[width=70pt]{ 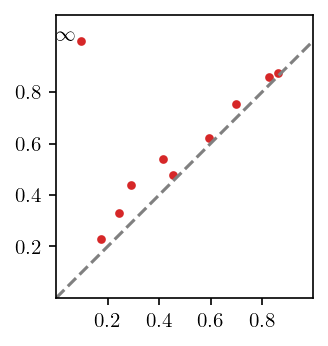}};
    \node at (1, -4.7) {\includegraphics[width=70pt]{ 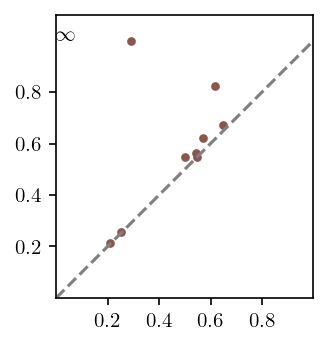}};
    \node at (3, -4.7) {\includegraphics[width=70pt]{ 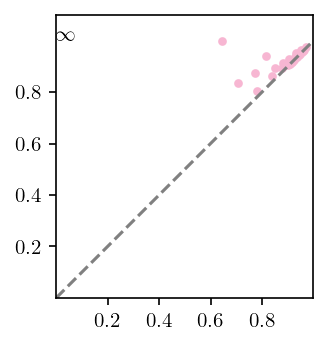}};
    \node at (0, -7.3) {
    \begin{tabular}{c|c|r|r|r|r}
        cluster & \#cities & ave. pop & ave. \# tracts & ave. \% HISP \\ 
1 \textcolor {mycolor1}{$\blacksquare$} & 4 & 4784129.25  & 1252.75  & 38.64 \% \\ 
2 \textcolor {mycolor2}{$\blacksquare$} & 11 & 1206203.55  & 301.09  & 28.50 \% \\ 
3 \textcolor {mycolor3}{$\blacksquare$} & 27 & 593933.30  & 145.81  & 49.63 \% \\ 
4 \textcolor {mycolor4}{$\blacksquare$} & 24 & 611312.54  & 158.75  & 18.23 \% \\ 
5 \textcolor {mycolor5}{$\blacksquare$} & 34 & 500748.00  & 138.26  & 9.18 \% \\ 
    \end{tabular}
    };
    \end{tikzpicture}
    \caption{One hundred U.S. cities clustered into $k = 5$ clusters based on their Hispanic \% persistence diagrams. Above, city locations are plotted and colored by cluster assignment. In the middle, the cluster mean for each cluster is shown. Each cluster mean is an average of the persistence diagrams for that cluster. At the bottom is a table collecting various statistics for the cities within each cluster.}
    \label{fig:HISP_map}
\end{figure}

The five clusters are shown in Figure \ref{fig:HISP_map}. Cluster 1 consists of large cities with a large number of clusters with high Hispanic populations e.g. New York NY, Los Angeles CA, and Houston TX. Cluster 2 also contains some large cities such as Chicago IL, but these cities typically have less high Hispanic percentage clusters than cities in cluster 1. Cluster 3 is similar to cluster 3 (for Black population) in Figure \ref{fig:BLACK_map}, in that these cities (e.g.~Albuquerque NM and Miami FL) typically have just one Hispanic cluster, which could indicate high segregation levels. Clusters 4 and 5 are cities with low Hispanic populations. 

We see a similar categorization of cities for both Black and Hispanic population, namely: (1) outliers with many clusters, (2) cities with a few clusters, (3) cities with only one cluster, (4) \& (5) cities with lower overall concentrations of Black/Hispanic population. In both maps we see some geographic clustering, i.e.~groups of nearby cities belonging to the same cluster. A lot of this effect is due to the fact that low Black/Hispanic population cities are clustered together into the last two clusters in each case, and large scale geography naturally tends to regulate the overall percentage of the population in these demographics. Nonetheless, it indicates that persistence data is able to pick up on some regional effects. We should also caution that the choice of $k$ and the initialization of the cluster means has a large effect on the final clustering, so that not too much should be read into the particular cluster labels produced in Figures \ref{fig:BLACK_map} and \ref{fig:HISP_map}. Rather, these labels can serve as jumping-off points for further analysis. For example, if nearby cities are in different clusters, then inspecting their persistence diagrams can confirm the difference, which can then lead to a detailed investigation of where this difference comes from.

\section{Modifiable Areal Unit Problem (MAUP)} \label{sec:drawback}

So far we have seen a number of instances in which our persistent homology method is able to flag conspicuous patterns in demographic data, which upon further investigation can be explained by a number of different factors such as changing demographics at both large and small scales, patterns of segregation, or even the location of correctional facilities. In each case, we have followed up the persistence diagram analysis with a direct inspection of the raw data. This is for good reason, as our method is also very capable of producing false positives due to anomalies in the data. The main source of this issue is the \emph{modifiable areal unit problem} (MAUP)\cite{Berry1968_Spatial}. Since percentages are based on Census tracts, the way tracts are drawn can have an effect on the resulting persistence diagram, even if the population remains completely static. This problem is not unique to our method, and shows up in almost any method where values are aggregated to specific geographic units. Measures of segregation and spatial clustering such as Moran's I~\cite{moran1950notes} are some notable examples of methods which are vulnerable to the MAUP. Attempts have been made to design methods which avoid the MAUP (see e.g.~\cite{hennerdal2017multiscalar}). However, for now we do not attempt to immunize our method against it, and instead take a brief look at some examples to demonstrate the kinds of MAUP effects that can be observed.\\

\subsection*{Example 1: Minneapolis MN} 
\begin{figure}
    \centering
    \begin{subfigure}{0.45\textwidth}
    \centering
    Minneapolis (2010)
            \includegraphics[width=\textwidth]{ 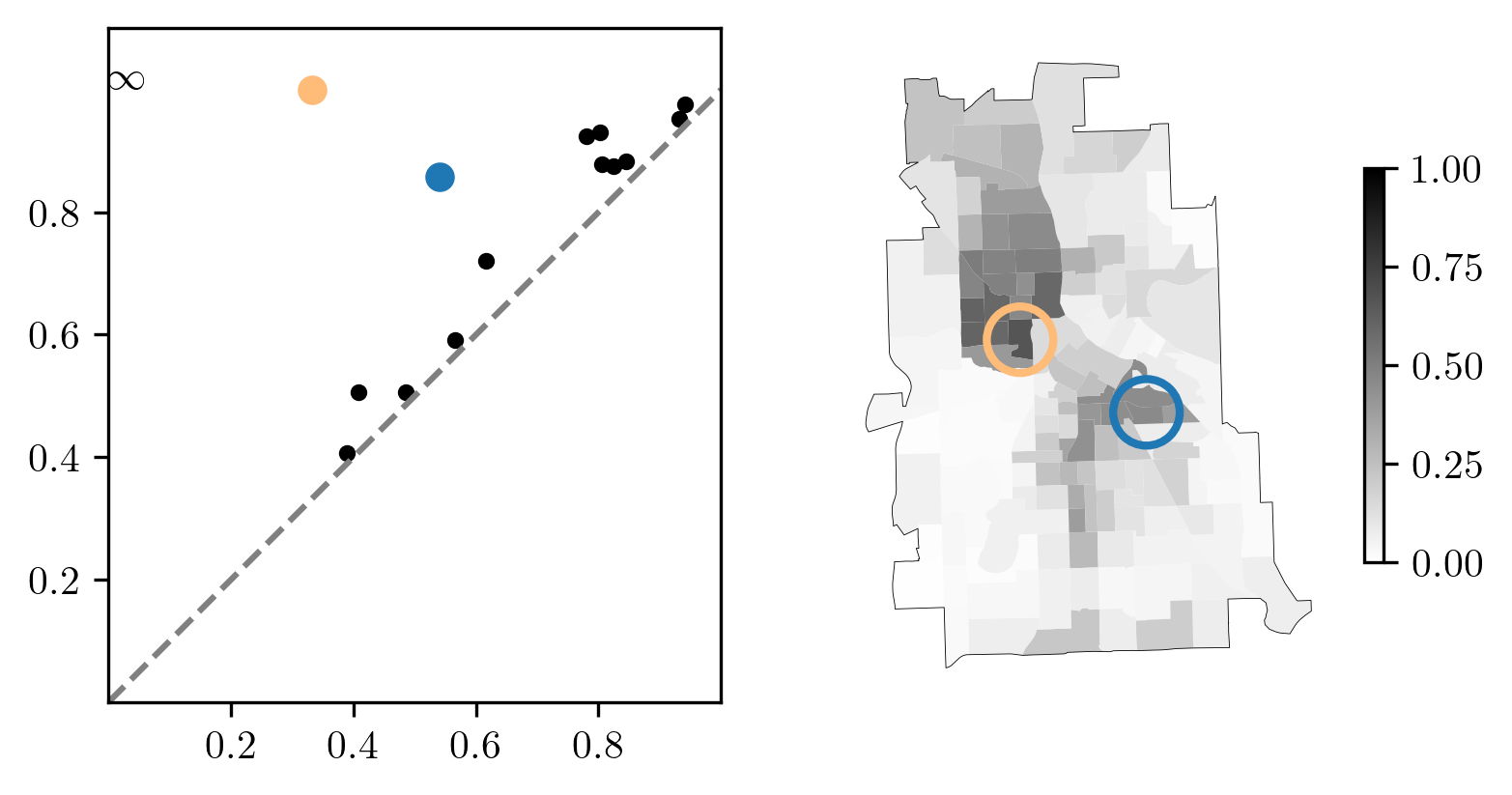}
    \end{subfigure}\hfill 
    \begin{subfigure}{0.45\textwidth}
    \centering
    Minneapolis (2020)
            \includegraphics[width=\textwidth]{ 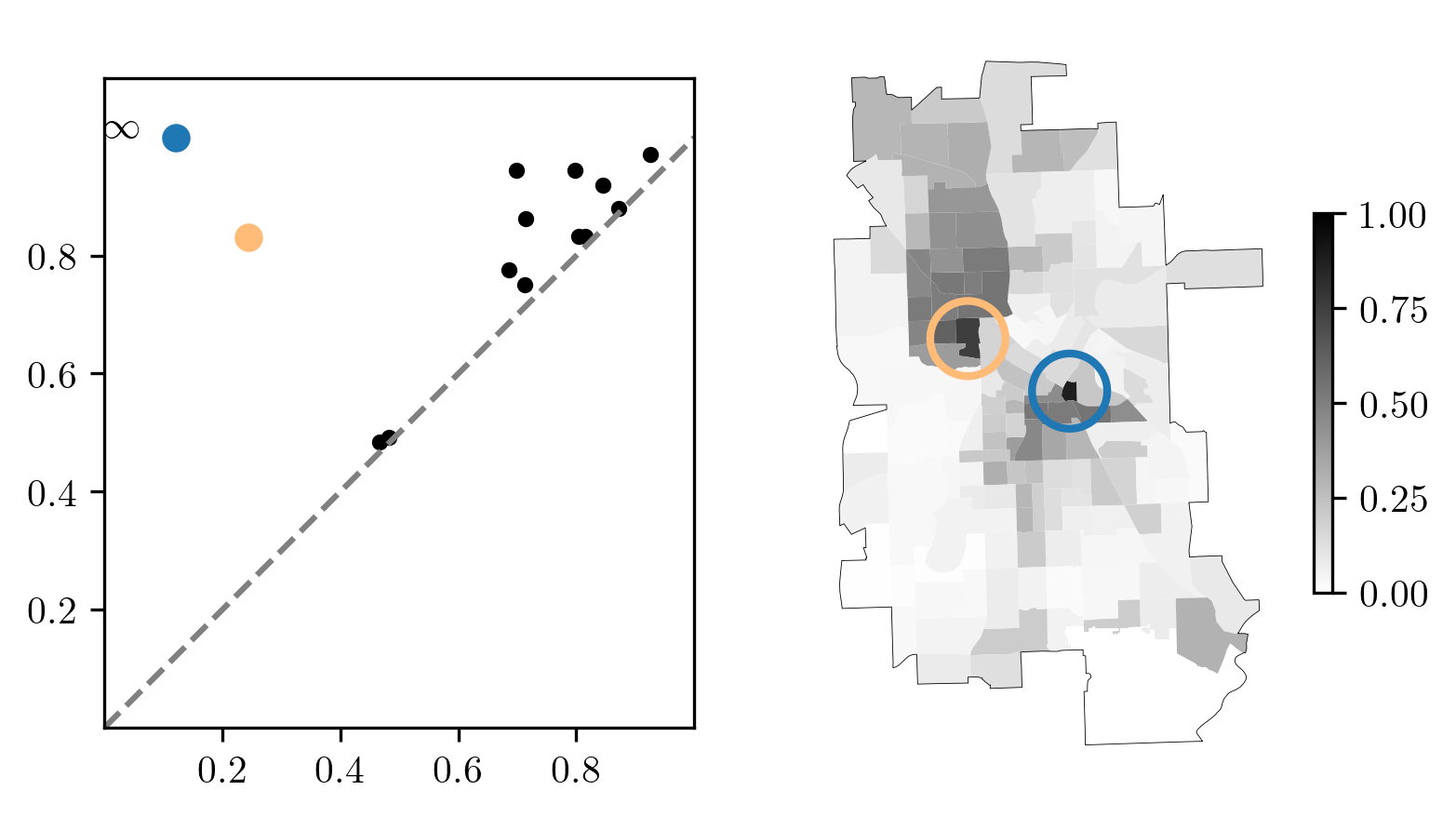}
    \end{subfigure}
    \caption{Contrasting Minneapolis MN in 2010 and in 2020 using Black population. One of the points (beige) moves due to a genuine population shift, the other (blue) is affected by a change in tracts between 2010 and 2020.}
    \label{fig:minneapolis}
\end{figure}

Figure \ref{fig:minneapolis} compares 2010 and 2020 Black population data for Minneapolis. Points are colored based on their birth tract location. The shift in the blue point's birth time is the result of a 45\% Black tract being broken up, with the area being included into a 88\% Black tract and a 22\% Black tract. Thus, this change is most likely an artifact of the way tracts were redrawn.

\subsection*{Example 2: San Diego} 
\begin{figure}
    \centering
    \begin{subfigure}{0.45\textwidth}
    \centering
    San Diego (2010)
            \includegraphics[width=\textwidth]{ 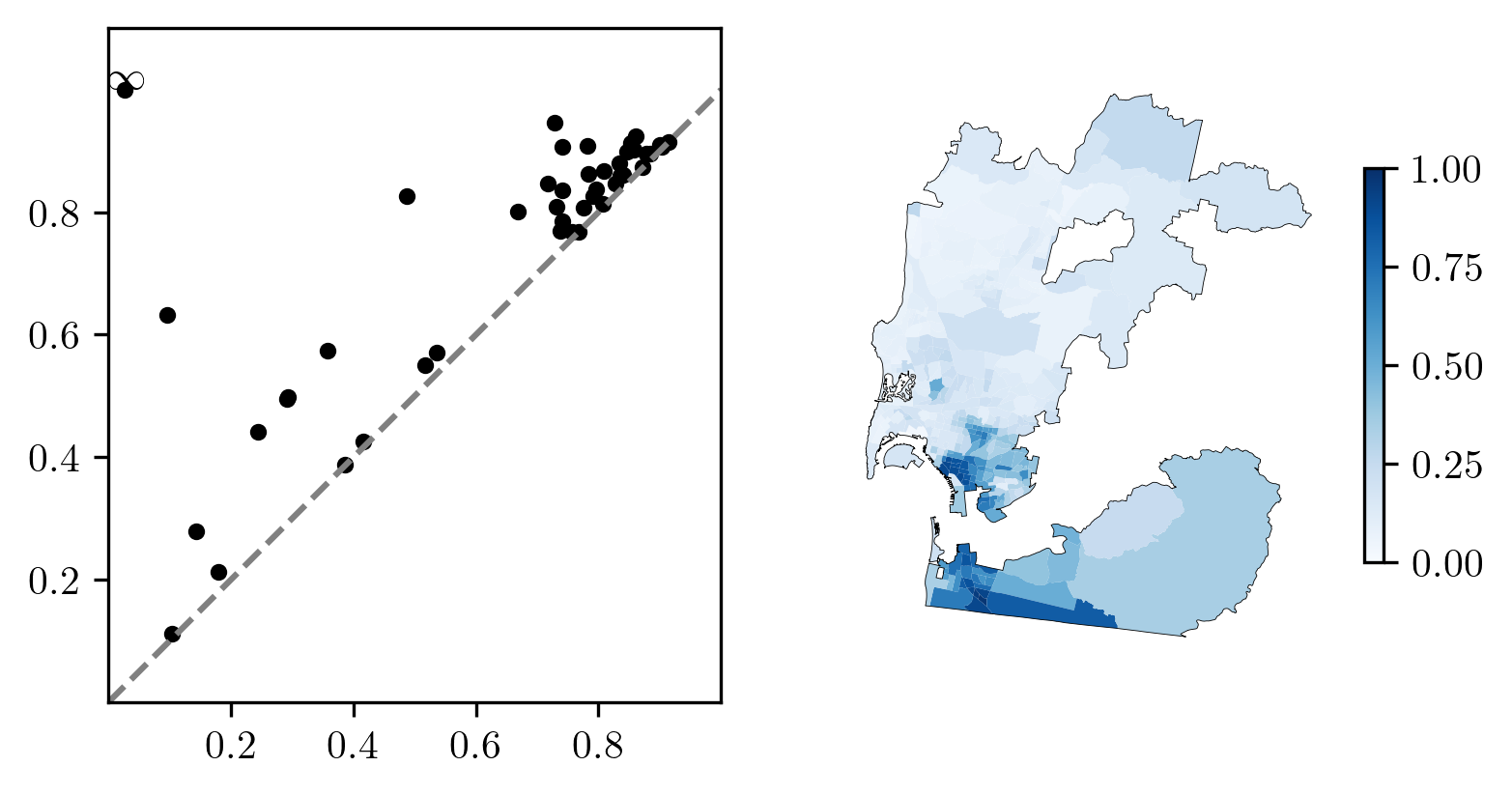}
    \end{subfigure}\hfill 
    \begin{subfigure}{0.45\textwidth}
    \centering
    San Diego (2020)
            \includegraphics[width=\textwidth]{ 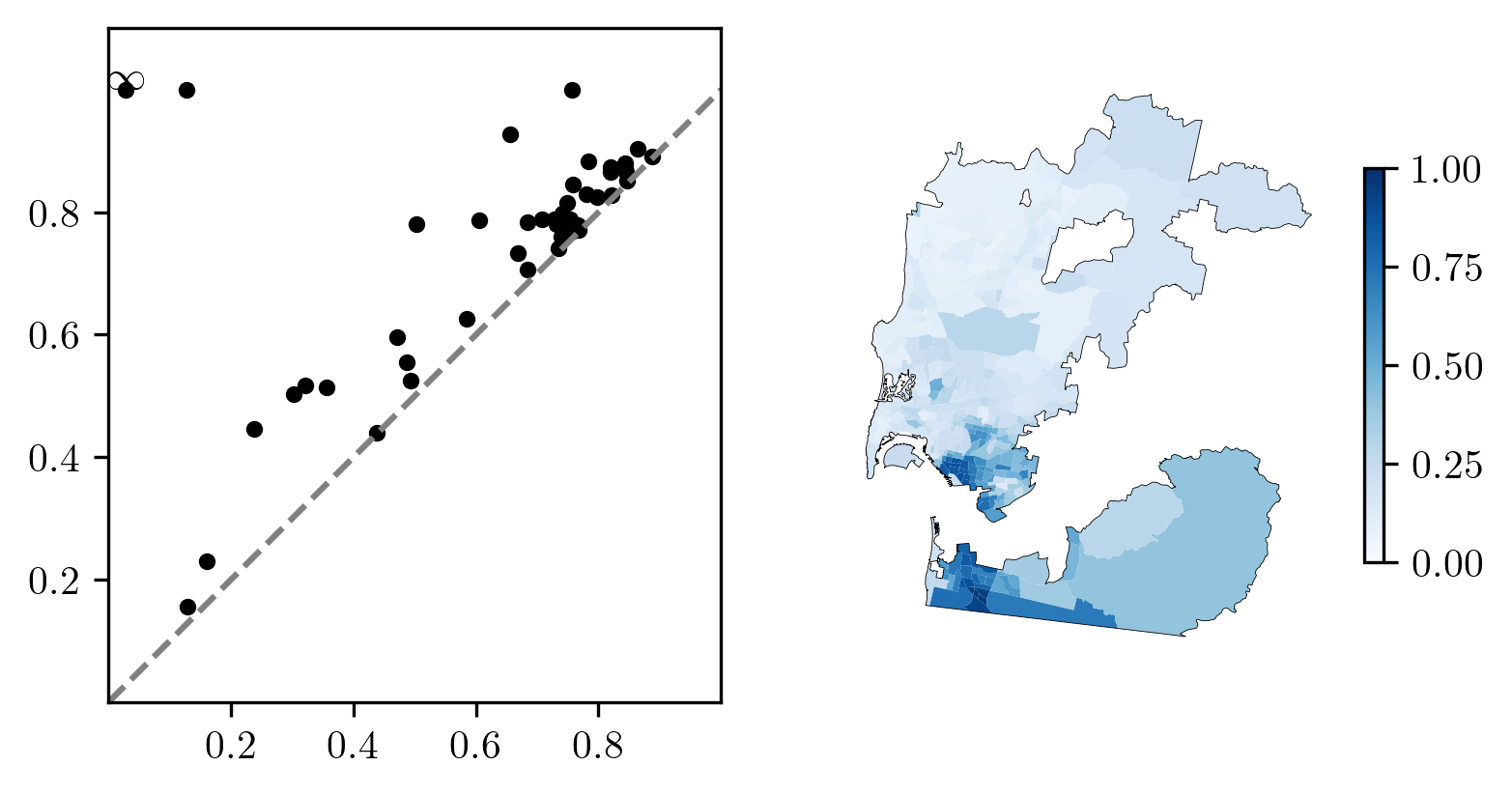}
    \end{subfigure}
    \caption{Contrasting San Diego CA in 2010 and in 2020 using Hispanic population. The addition of a point at infinity is due to the dual graph being connected in 2010 but disconnected in 2020.}
    \label{fig:sandiego}
\end{figure}

From 2010 to 2020, the persistence diagram for Hispanic population in San Diego gained an extra point at infinity, as shown in Figure~\ref{fig:sandiego}. The reason (not visible on the maps) is that in 2010, a tract which reached across the San Diego Bay was broken up so that the dual graph went from being connected in 2010 to disconnected in 2020. There is no well-defined way to build dual graphs in the presence of water features or islands, something that causes serious problems for redistricting simulations~\cite{caldera2020mathematics}. We chose not to fill in edges by hand across water features in our analysis. 

\section{Conclusion}
Our goal has been to showcase some of the advantages and drawbacks of persistent homology as an exploratory analysis tool for studying demographic data, adding to the growing literature on applications of persistent homology to (geo)spatial data~\cite{feng2020spatial}. We have demonstrated how persistent homology provides a range of quantitative surveys of 100 cities and their demographic patterns. Nonetheless, there remain some obstacles and open questions which must be resolved in order to increase the value of persistent homology methods for this kind of data. These include mitigating the effect of the MAUP, dealing with separation by water, more robust clustering techniques, and investigating other persistence statistics beyond total persistence. Our work also underlines the need for following up persistent homology analysis with a closer look at the underlying data, guided by the location information connected to the persistence diagram. Persistence diagrams, on their own, can obscure important information or even be misleading.

Perhaps the most important objective not addressed in this paper is to connect the quantitative analysis provided by persistent homology to concepts from the social sciences. For example, persistent homology is connected to the fifth dimension (clustering) of segregation as defined by Massey and Denton~\cite{massey1988dimensions}, but it is also capable of detecting much more. Detecting demographic shifts with persistent homology is just the first step, and should lead to a search for the social forces behind these shifts, relying on the established ideas and methods outside of mathematics.

\bibliographystyle{plain}
\bibliography{myreferences}

\end{document}